\documentclass[12pt]{article}
\usepackage{amsfonts, amssymb, amsmath, amsthm}

\makeatletter
\newtheorem{theorem}{Theorem}
\@addtoreset{theorem}{section}

\newtheorem{lemma}{Lemma}
\@addtoreset{lemma}{section}

\@addtoreset{definition}{section}

\newtheorem{remark}{Remark}
\@addtoreset{remark}{section}

\newtheorem{example}{Example}
\@addtoreset{example}{section}

\@addtoreset{condition}{section}

\@addtoreset{corollary}{section}

\@addtoreset{equation}{section}

\renewcommand{\Im}{{\rm Im\,}}

\renewcommand{\ker}{{\rm ker\,}}
\renewcommand{\dim}{{\rm dim\,}}
\newcommand{\codim}{{\rm codim\,}}
\newcommand{\arctg}{{\rm arctg\,}}

\newcommand{\rank}{{\rm rank\,}}

\renewcommand{\@makefnmark}{\hbox{\mathsurround=0pt $^{\@thefnmark)}$}}
\renewcommand{\@makefntext}[1]{\parindent=1em\noindent
           \hbox to 1.8em{\hss$^{\@thefnmark)}$}#1}

\begin{document}
\setcounter{section}{0}
\setcounter{footnote}{0}

%

%

\begin{center}
{\Large Asymptotics of Solutions for Nonlocal Elliptic Problems\\
in Plane Angles} \footnote{This work was partially supported by
Russian Foundation for Basic Research (grant No~01-01-01030) and
by Russian Ministry for Education (grant No~E00-1.0-195).}

Pavel Gurevich
\end{center}
\renewcommand{\section}{\@startsection{section}{2}{0pt}{-3.5ex plus -1ex minus -.2ex}{1ex}{\bf}}

\begin{abstract}
We investigate asymptotic behavior of solutions for nonlocal
elliptic boundary value problems in plane angles and in~${\mathbb
R}^2\backslash\{0\}$. Such problems arise as model ones when
studying asymptotics of solutions for nonlocal elliptic problems
in bounded domains. We obtain explicit formulas for the asymptotic
coefficients in terms of eigenvectors and associated vectors of
both adjoint nonlocal operators (acting in spaces of
distributions) and formally adjoint (with respect to the Green
formula) nonlocal problems.
\end{abstract}

\section{Introduction}
T.~Carleman~\cite{Carleman} was one of the first who began studying nonlocal elliptic problems.
Investigation of nonlocal problems with shifts mapping a boundary on itself are
closely associated with paper~\cite{Carleman}. In~\cite{BitzSam}, A.V.~Bitsadze and A.A.~Samarskii
considered the Laplace equation in a domain $G\subset{\mathbb R}^n$ with the boundary condition that
connects the values of an unknown function on a manifold $\Upsilon_1\subset\partial G$ with its values 
on some manifold inside $G$; on $\partial G\backslash \Upsilon_1$ the Dirichlet condition was set.
Such a formulation is associated with further investigating nonlocal problems with shifts mapping a boundary
inside a domain. One can find a detailed bibliography devoted to nonlocal elliptic problems in~\cite{SkBook}.

In the theory of nonlocal elliptic problems of this type, the most difficult case
deals with the situation when support of nonlocal terms intersects with a boundary~\cite{BitzDAN85}--\cite{GM}.
This leads to appearance of power singularities for solutions near some set~${\cal K}$. Therefore, 
the problem of asymptotics of solutions near this set arises. Asymptotic formulas for solutions to
nonlocal elliptic problems in plane domains were first obtained by A.L.~Skubachevskii in~\cite{SkMs86}. 
They allow to prove a number of principally new properties (in comparison with ``local" elliptic problems both
in domains with angular points~\cite{KondrTMMO67,KondrOleinik} and in domains with smooth boundary). 
For example, smoothness of generalized solutions for nonlocal elliptic problems can be violated both near
vertexes of  small angles and near smooth boundary even for arbitrary small coefficients in nonlocal 
terms~\cite{SkMs86, SkRJMP}.

In this paper we investigate asymptotic behavior of solutions for nonlocal elliptic boundary value problems 
in plane angles and in~${\mathbb R}^2\backslash\{0\}$. Such problems arise as model ones when
studying asymptotics of solutions for nonlocal elliptic problems in bounded domains near the set~${\cal K}$.
We obtain explicit formulas for the 
asymptotic coefficients in terms of eigenvectors and associated vectors of both adjoint nonlocal operators
(acting in spaces of distributions) and formally adjoint (with respect to the Green formula) nonlocal problems.
Earlier adjoint nonlocal problems were studied in~\cite{GurGiess, GurDAN}.

Notice that a number of statements are proved similarly to results of papers~\cite{NP, MP}. In these cases
we shall give just schemes of proofs.
\section{Statement of nonlocal problems in plane angles and preliminary information. Asymptotics of solutions}\label{sectStatement}

{\bf 1.} Consider the plane angle $K=\{y\in{\mathbb R}^2:\ r>0,\ b_1<\omega<b_2\}$ with the sides
$\gamma_\sigma=\{y\in{\mathbb R}^2:\ r>0,\ \omega=b_\sigma\}\ (\sigma=1,\ 2)$.
Here $(\omega,\ r)$ are polar coordinates of a point $y;$ $-\pi<b_1<b_2<\pi.$

Denote by ${\cal P}(D_y),$ $B_{\sigma\mu}(D_y)$ and $B^{\cal G}_{\sigma\mu}(D_y)$
homogeneous differential operators with constant complex coefficients of orders
$2m,$  $m_{\sigma\mu}\le 2m-1$, and
$m_{\sigma\mu}\le 2m-1$ correspondingly ($\sigma=1,\ 2;$ $\mu=1,\ \dots,\ m$).
We shall suppose that the operator ${\cal P}(D_y)$ is properly elliptic and the system of operators
$\{B_{\sigma\mu}(D_y)\}_{\mu=1}^m$ is normal and covers ${\cal P}(D_y)$ on $\gamma_\sigma$ (see~\cite[Chapter 2]{LM}).
We do not impose any conditions (except the restrictions on orders) on the operators $B^{\cal G}_{\sigma\mu}(D_y)$,
which play further the role of nonlocal ones.

Consider the following nonlocal elliptic problem in the plane angle~$K$:
\begin{equation}\label{eqP}
  {\cal P}(D_y)u=f(y) \quad (y\in K),
\end{equation}
\begin{equation}\label{eqB}
  \begin{array}{c}
  {\cal B}_{\sigma\mu}(D_y)u\equiv B_{\sigma\mu}(D_y)u(y)|_{\gamma_\sigma}+
                   (B^{\cal G}_{\sigma\mu}(D_y)u)({\cal G}_\sigma y)|_{\gamma_\sigma}
    =g_{\sigma\mu}(y) \quad (y\in\gamma_\sigma),\\
  \sigma=1,\ 2;\ \mu=1,\ \dots,\ m.
  \end{array}
\end{equation}
The notation $(B^{\cal G}_{\sigma\mu}(D_y)u)({\cal G}_\sigma y)$ means that the expression
$(B^{\cal G}_{\sigma\mu}(D_{y'})u)(y')$ is taken for $y'={\cal G}_\sigma y;$
${\cal G}_\sigma$ is the operator of rotation by the angle~$\omega_\sigma$ and expansion by~$\beta_\sigma$
times in the plane~$\{y\}$ such that
$b_1<b_1+\omega_1=b_2+\omega_2=b<b_2,$ $0<\beta_\sigma.$

For any set $G\subset{\mathbb R}^n$ ($n\ge1$), denote by~$C_0^\infty(G)$ the set of infinitely differentiable
in~$\bar G$ functions with supports belonging to~$G$.
We introduce the space $H_a^l(K)$ as a completion of the set $C_0^\infty(\bar K\backslash \{0\})$
in the norm
$
 \|w\|_{H_a^l(K)}=\Bigl(
    \sum\limits_{|\alpha|\le l}\int\limits_K r^{2(a-l+|\alpha|)} |D_y^\alpha w(y)|^2 dy
                                       \Bigr)^{1/2},
$
where $a\in\mathbb R,$ $l\ge 0$ is an integer. By $H_a^{l-1/2}(\gamma')$ for $l\ge1$ we denote the space of traces on a ray
$\gamma'=\{y\in{\mathbb R}^2:\ r>0,\ \omega=b'\}\ (b_1\le b'\le b_2)$ with the norm
$
 \|\psi\|_{H_a^{l-1/2}(\gamma')}=\inf\|w\|_{H_a^l(K)}\ (w\in H_a^l(K):\  w|_{\gamma'} = \psi).
$

Introduce the bounded operator corresponding to problem~(\ref{eqP}), (\ref{eqB})
$$
 \begin{array}{c}
 {\cal L}=\{{\cal P}(D_y),\ {\cal B}_{\sigma\mu}(D_y)\}: H_{a}^{l+2m}(K)\to \\
 \to H_a^l(K,\ \gamma)=
 H_a^l(K)\times\prod\limits_{\sigma=1,2}\prod\limits_{\mu=1}^m H_{a}^{l+2m-m_{\sigma\mu}-1/2}(\gamma_\sigma).
 \end{array}
$$

{\bf 2.}
Write the operators ${\cal P}(D_y),$ $B_{\sigma\mu}(D_y),$ $B^{\cal G}_{\sigma\mu}(D_y)$ in polar coordinates:
${\cal P}(D_y)=r^{-2m}\tilde{\cal P}(\omega,\ D_\omega,\ rD_r),$
$B_{\sigma\mu}(D_y)=r^{-m_{\sigma\mu}} \tilde B_{\sigma\mu}(\omega,\ D_\omega,\ rD_r),$
$B^{\cal G}_{\sigma\mu}(D_y)=r^{-m_{\sigma\mu}}\tilde B^{\cal G}_{\sigma\mu}(\omega,\ D_\omega,\ rD_r),$ where
$D_\omega=-i\frac{\displaystyle\partial}{\displaystyle\partial\omega},$
$D_r=-i\frac{\displaystyle\partial}{\displaystyle\partial r}.$

We shall denote by $\tilde w(\lambda)$ the Mellin transformation of a function $w\in C_0^\infty({\mathbb R}_+)$:
$$
 \tilde w(\lambda)=(2\pi)^{-1/2}\int\limits_0^\infty r^{-i\lambda-1}w(r)\,dr.
$$

Put $\{f,\ g_{\sigma\mu}\}=0$ in~(\ref{eqP}) and (\ref{eqB}) and do formally the Mellin transformation. Then we get
\begin{equation}\label{eqPLambda}
  \tilde{\cal P}(\lambda) \tilde u(\omega,\ \lambda)=0 \quad
   (b_1<\omega<b_2),
\end{equation}
\begin{equation}\label{eqBLambda}
 \begin{array}{c}
   \tilde {\cal B}_{\sigma\mu}(\lambda) \tilde u(\omega,\ \lambda)\equiv
  \tilde B_{\sigma\mu}(\lambda) \tilde u(\omega,\ \lambda)|_{\omega=b_\sigma}
              +\beta_\sigma^{-m_{\sigma\mu}+i\lambda} \tilde B^{\cal G}_{\sigma\mu}(\lambda)
              \tilde u(\omega+\omega_\sigma,\ \lambda)|_{\omega=b_\sigma}=0.
  \end{array}
\end{equation}
Here (and further) we omit for short the arguments~$\omega$ and $D_\omega$ in differential operators.
This problem is ordinary differential equation~(\ref{eqPLambda})
with nonlocal conditions~(\ref{eqBLambda}) that connect the values of a solution~$\tilde u$ and its derivatives
at the point~$\omega=b_\sigma$ with the values of a solution~$\tilde u$ and its derivatives at the internal
point~$\omega=b$ of the interval~$(b_1,\ b_2).$ Asymptotics of solutions for nonlocal
problem~(\ref{eqP}), (\ref{eqB}) in the angle~$K$ will be described in terms of eigenvalues and corresponding Jordan chains
of problem~(\ref{eqPLambda}), (\ref{eqBLambda}).

Let us consider the operator--valued function corresponding to nonlocal problem~(\ref{eqPLambda}), (\ref{eqBLambda})
$$
 \begin{array}{c}
 \tilde{\cal L}(\lambda)=\{\tilde{\cal P}(\lambda),\
  \tilde{\cal B}_{\sigma\mu}(\lambda)\}:W^{l+2m}(b_1,\ b_2)
  \to W^{l}[b_1,\ b_2]=W^l(b_{1},\ b_{2}) \times{\mathbb C}^{2m},
 \end{array}
$$
Here $W^l(\cdot)=W^l_2(\cdot)$ is the Sobolev space of order $l\ge 0$ (if $l=0$, we put $W^0(\cdot)=L_2(\cdot)$).

Now we shall remind some well-known definitions and facts
(see~\cite{GS}). A holomorphic at a point $\lambda_0$
vector--function~$\varphi(\lambda)$ with the values
in~$W^{l+2m}(b_1,\ b_2)$ is called {\it a root function} of the
operator~$\tilde{\cal L}(\lambda)$ at $\lambda_0$ if
$\varphi(\lambda_0)\ne0$ and the vector--function~$\tilde{\cal
L}(\lambda)\varphi(\lambda)$ is equal to zero at~$\lambda_0$.
If~$\tilde{\cal L}(\lambda)$ has at least one root function at a
point~$\lambda_0$, then $\lambda_0$ is called {\it an eigenvalue}
of $\tilde{\cal L}(\lambda)$. Multiplicity of zero for the
vector--function~$\tilde{\cal L}(\lambda)\varphi(\lambda)$ at the
point~$\lambda_0$ as called {\it a multiplicity of the root
function} $\varphi(\lambda)$; the
vector~$\varphi^{(0)}=\varphi(\lambda_0)$ is called {\it an
eigenvector} corresponding to the eigenvalue~$\lambda_0$. Let
$\varphi(\lambda)$ be a root function at a point~$\lambda_0$ of
multiplicity~$\varkappa$ and
$\varphi(\lambda)=\sum\limits_{j=0}^{\infty}(\lambda-\lambda_0)^j\varphi^{(j)}$.
Then the vectors $\varphi^{(1)},\ \dots,\ \varphi^{(\varkappa-1)}$
are called {\it associated with the eigenvector} $\varphi_0$, and
the ordered set~$\varphi^{(0)},\ \dots,\ \varphi^{(\varkappa-1)}$
is called  {\it a Jordan chain} corresponding to the
eigenvalue~$\lambda_0$. {\it Rank} of the
eigenvector~$\varphi^{(0)}$ ($\rank\varphi^{(0)}$) is the maximum
of multiplicities of all root functions such that
$\varphi(\lambda_0)=\varphi^{(0)}$.

\begin{remark}\label{rSmoothEigVecL}
An eigenvector and associated vectors $\varphi^{(0)},\ \dots,\ \varphi^{(\varkappa-1)}$ of the operator
$\tilde{\cal L}(\lambda)$ corresponding to an eigenvalue~$\lambda_0$ satisfy the equalities
\begin{equation}\label{eqDefEigVecL}
 \sum\limits_{q=0}^\nu\frac{\displaystyle 1}{\displaystyle q!}\partial_\lambda^q \tilde{\cal L}(\lambda_0)\varphi^{(\nu-q)}=0,\quad
 \nu=0,\ \dots,\ \varkappa-1.
\end{equation}
Here and further $\partial_\lambda^q$ is the derivative of order~$q$ with respect to $\lambda$.

From equalities~(\ref{eqDefEigVecL}) and Lemma~\ref{lSmoothLLambda}, it follows that
eigenvectors and associated vectors of the operator~$\tilde{\cal L}(\lambda)$ are infinitely differentiable functions
in the interval~$[b_1,\ b_2]$.
\end{remark}

From Lemma~2.1~\cite{SkDu90}, it follows that all eigenvalues of the operator~$\tilde{\cal L}(\lambda)$ are isolated.
Moreover, $\dim\ker\tilde{\cal L}(\lambda_0)<\infty$ for any eigenvalue~$\lambda_0$, and ranks of all eigenvectors are finite.
Suppose $J=\dim\ker\tilde{\cal L}(\lambda_0)$ and $\varphi^{(0,1)},\ \dots,\ \varphi^{(0,J)}$ is a system of
eigenvectors such that $\rank\varphi^{(0,1)}$ is the greatest of ranks of all eigenvectors corresponding to the
eigenvalue~$\lambda_0$, and $\rank\varphi^{(0,j)}$ ($j=2,\ \dots,\ J$) is the greatest of ranks of eigenvectors from
some orthogonal supplement in~$\ker\tilde{\cal L}(\lambda_0)$ to the linear manifold of the
vectors~$\varphi^{(0,1)},\ \dots,\ \varphi^{(0,j-1)}$. The numbers~$\varkappa_j=\rank\varphi^{(0,j)}$ are called
{\it partial multiplicities} of the eigenvalue~$\lambda_0$, and the sum~$\varkappa_1+\dots+\varkappa_J$ is called
{\it a (full) multiplicity} of $\lambda_0$.
If the vectors~$\varphi^{(0,j)},\ \dots,\ \varphi^{(\varkappa_j-1,j)}$ form a Jordan chain for every~$j=1,\ \dots,\ J$,
then the set of vectors $\{\varphi^{(0,j)},\ \dots,\ \varphi^{(\varkappa_j-1,j)}:j=1,\ \dots,\ J\}$ is called
{\it a canonical system of Jordan chains} corresponding to the eigenvalue~$\lambda_0$.

\begin{example}\label{exEV}{\rm
Put $b_1=-\omega_0,\ b_2=\omega_0$.
In the plane angle~$K=\{y\in{\mathbb R}^2: |\omega|<\omega_0\}$ ($0<\omega_0<\pi$) with the sides
$\gamma_\sigma=\{y\in{\mathbb R}^2: \omega=(-1)^\sigma\omega_0\}$, $\sigma=1,\ 2$, we consider the
nonlocal problem
\begin{equation}\label{eqPEV}
  \triangle u=f(y)\quad (y\in K),
\end{equation}
\begin{equation}\label{eqBEV}
  u|_{\gamma_1}=0,\   u|_{\gamma_2}+b\,  u({\cal G}_2y)|_{\gamma_2}=0,
\end{equation}
where $b\in{\mathbb R}$, ${\cal G}_2$ is the operator of rotation by the angle~$-\omega_0$. The following
model nonlocal eigenvalue problem corresponds to problem~(\ref{eqPEV}), (\ref{eqBEV}):
\begin{equation}\label{eqPEVLambda}
  \frac{d^2\varphi(\omega)}{d\omega^2}-\lambda^2\varphi(\omega)=0\quad (|\omega|<\omega_0),
\end{equation}
\begin{equation}\label{eqBEVLambda}
  \varphi(-\omega_0)=0,\   \varphi(\omega_0)+b\, \varphi(0)=0.
\end{equation}

One can immediately check (see also~\cite[Chapter 2]{NP}) that, for~$b=0$, (that is if
problem~(\ref{eqPEV}), (\ref{eqBEV}) is ``local") the eigenvalues of
problem~(\ref{eqPEVLambda}), (\ref{eqBEVLambda}) have the form~$\lambda_k=i\frac{\pi k}{2\omega_0}$,
$k\in{\mathbb Z}\backslash\{0\}$;
The eigenvectors $\varphi_k^{(0)}(\omega)=e^{i\frac{\pi k}{2\omega_0}\omega}-e^{-i\frac{\pi k}{2\omega_0}\omega}$
correspond to these eigenvalues. Associated vectors are absent, that is all the eigenvalues are of multiplicities 1.

Now we shall show that, for $b\ne0$, there may be Jordan chains with a length more than 1 corresponding
to eigenvalues of problem~(\ref{eqPEVLambda}), (\ref{eqBEVLambda}).

I) First we consider the case~$\lambda\ne 0$. Substituting the general solution
$\varphi(\omega)=c_1e^{\lambda\omega}+c_2e^{-\lambda\omega}$ for equation~(\ref{eqPEVLambda}) into
nonlocal conditions~(\ref{eqBEVLambda}), we get
\begin{equation}\label{eqSys0}
\begin{array}{ll}
 c_1e^{-\lambda\omega_0}+c_2e^{\lambda\omega_0}=0,\\
 (e^{\lambda\omega_0}+b)c_1+(e^{-\lambda\omega_0}+b)c_2=0.
\end{array}
\end{equation}
Equate the determinant~$D(\lambda)$ of system~(\ref{eqSys0}) with zero:
$$
 (e^{-\lambda\omega_0}-e^{\lambda\omega_0})(e^{\lambda\omega_0}+e^{-\lambda\omega_0}+b)=0.
$$

1) Let we have $e^{-\lambda\omega_0}-e^{\lambda\omega_0}=0$. Then we obtain the series of eigenvalues
$$
 \lambda_{1k}=i\frac{\pi k}{\omega_0},\ k\in{\mathbb Z}\backslash\{0\};
$$
The eigenvectors
$$
 \varphi_{1k}^{(0)}(\omega)=e^{i\frac{\pi k}{\omega_0}\omega}-e^{-i\frac{\pi k}{\omega_0}\omega}
$$
correspond to these eigenvalues.
Consider a problem of finding an associated vector~$\varphi_{1k}^{(1)}(\omega)$.
According to~(\ref{eqDefEigVecL}),  $\varphi_{1k}^{(1)}(\omega)$ satisfies the equation
$$
 \frac{d^2\varphi_{1k}^{(1)}}{d\omega^2}+\frac{(\pi k)^2}{\omega_0^2}\varphi_{1k}^{(1)}
 -2i\frac{\pi k}{\omega_0}\varphi_{1k}^{(0)}=0\quad (|\omega|<\omega_0)
$$
and nonlocal conditions~(\ref{eqBEVLambda}). Substituting the general solution
$\varphi(\omega)=c_1e^{i\frac{\pi k}{\omega_0}\omega}+c_2e^{-i\frac{\pi k}{\omega_0}\omega}+
\omega(e^{i\frac{\pi k}{\omega_0}\omega}+e^{-i\frac{\pi k}{\omega_0}\omega})$ of the last equation into
nonlocal conditions~(\ref{eqBEVLambda}), we get
\begin{equation}\label{eqSys1k}
\begin{array}{ll}
 c_1+c_2=2\omega_0,\\
 ((-1)^k+b)c_1+((-1)^k+b)c_2=-2(-1)^k\omega_0.
\end{array}
\end{equation}
Therefore, an associated vector~$\varphi_{1k}^{(1)}(\omega)$ exists if and only if
$$
 b=2(-1)^{k+1}.
$$
If $b=2(-1)^{k+1}$, we can put
$$
 \varphi_{1k}^{(1)}(\omega)=(\omega+2\omega_0)e^{i\frac{\pi k}{\omega_0}\omega}+\omega e^{-i\frac{\pi k}{\omega_0}\omega}.
$$
Analogously, using~(\ref{eqDefEigVecL}), we find the second associated vector
$$
  \varphi_{1k}^{(2)}(\omega)=(\frac{\omega^2}{2}+2\omega_0\omega+2\omega_0^2)e^{i\frac{\pi k}{\omega_0}\omega}
  -\frac{\omega^2}{2} e^{-i\frac{\pi k}{\omega_0}\omega}.
$$
One can directly check that the third associated vector is absent.

2) Let we have
\begin{equation}\label{eqCase2}
e^{\lambda\omega_0}+e^{-\lambda\omega_0}+b=0.
\end{equation}
Then we obtain the following series of eigenvalues:
$$
\begin{array}{c}
 \lambda_{2n}^\pm=\frac{\ln\left(-\frac{b}{2}\pm\frac{\sqrt{b^2-4} }{2}\right)}{\omega_0} +i\frac{2\pi n}{\omega_0}\ \mbox{for } b<-2;\
 \lambda_{2n}^\pm=i\frac{\pm\arctg\frac{\sqrt{4-b^2} }{b}+2\pi n}{\omega_0}\ \mbox{for } -2<b<0;\\
 \lambda_{2n}^\pm=i\frac{\pm\arctg\frac{\sqrt{4-b^2} }{b}+(2n+1)\pi}{\omega_0}\ \mbox{for } 0<b<2;\
 \lambda_{2n}^\pm=\frac{\ln\left(\frac{b}{2}\pm\frac{\sqrt{b^2-4} }{2}\right)}{\omega_0} +i\frac{(2n+1)\pi}{\omega_0}\
  \mbox{for } b>2;
\end{array}
$$
$n\in{\mathbb Z}$. If $|b|=2$, then we have eigenvalues from the
series~$\{\lambda_{1k}\}_{k\in{\mathbb Z}\backslash\{0\}}$, which is considered above.
The eigenvector
$$
 \varphi_{2n}^{(0)\pm}(\omega)=e^{\lambda_{2n}^\pm\omega}-e^{-2\lambda_{2n}^\pm\omega_0}e^{-\lambda_{2n}^\pm\omega}.
$$
corresponds to the eigenvalue~$\lambda_{2n}^\pm$. Let us show that there are no associated vectors
if~$\lambda=\lambda_{2n}^\pm$. Substitute the general solution
 $\varphi_{2n}^{(1)\pm}(\omega)=c_1e^{\lambda_{2n}^\pm\omega}+c_2e^{-\lambda_{2n}^\pm\omega}+
\omega(e^{\lambda_{2n}^\pm\omega}+e^{-2\lambda_{2n}^\pm\omega_0}e^{-\lambda_{2n}^\pm\omega})$ for the equation
$$
 \frac{d^2\varphi_{2n}^{(1)\pm}}{d\omega^2}-(\lambda_{2n}^\pm)^2 \varphi_{2n}^{(1)\pm}
 -2\lambda_{2n}^\pm \varphi_{2n}^{(0)\pm}=0\quad (|\omega|<\omega_0)
$$
into nonlocal conditions~(\ref{eqBEVLambda}). Then we have
\begin{equation}\label{eqSys2n}
\begin{array}{ll}
 e^{-\lambda_{2n}^\pm\omega_0}c_1+e^{\lambda_{2n}^\pm\omega_0}c_2=2\omega_0e^{-\lambda_{2n}^\pm\omega_0},\\
 (e^{\lambda_{2n}^\pm\omega_0}+b)c_1+(e^{\lambda_{2n}^\pm\omega_0}+b)c_2=
 -\omega_0(e^{\lambda_{2n}^\pm\omega_0}+e^{-3\lambda_{2n}^\pm\omega_0}).
\end{array}
\end{equation}
Rank of the matrix of system~(\ref{eqSys2n}) is equal to~1. Therefore, system~(\ref{eqSys2n}) is compatible
if and only if
$$
\left|
\begin{array}{cc}
 e^{-\lambda_{2n}^\pm\omega_0} & 2\omega_0e^{-\lambda_{2n}^\pm\omega_0}\\
 e^{\lambda_{2n}^\pm\omega_0}+b & -\omega_0(e^{\lambda_{2n}^\pm\omega_0}+e^{-3\lambda_{2n}^\pm\omega_0})
\end{array}
\right|=0.
$$
The last equality is equivalent to the following one:
$$
  3e^{\lambda_{2n}^\pm\omega_0}+e^{-3\lambda_{2n}^\pm\omega_0}+2b=0.
$$
From this, taking into account~(\ref{eqCase2}), it follows that either $e^{\lambda_{2n}^\pm\omega_0}=1$, $b=-2$ or
$e^{\lambda_{2n}^\pm\omega_0}=-1$, $b=2$. But we now consider the case~$|b|\ne2$. Hence, there are no
associated vectors if~$\lambda=\lambda_{2n}^\pm$.

II) The case~$\lambda=0$ is studied analogously. It turns out that $\lambda=0$ is an eigenvalue of
problem~(\ref{eqPEVLambda}), (\ref{eqBEVLambda}) if and only if~$b=-2$. Moreover, if $b=-2$, then for
the eigenvalue~$\lambda=0$, there exist one
eigenvector~$\varphi_0^{(0)}(\omega)=\omega+\omega_0$ and one associated vector~$\varphi_0^{(1)}(\omega)=0$.

Thus, we have shown that problem~(\ref{eqPEVLambda}), (\ref{eqBEVLambda}) has eigenvalues of multiplicities
 more than~1 if and only if~$|b|=2$.
}
\end{example}

{\bf 3.} The following result on isomorphism follows from~\cite[\S2]{SkDu90}.
\begin{theorem}\label{thSolv}
 Suppose the line~$\Im\,\lambda=a+1-l-2m$ contains no eigenvalues of the
 operator~$\tilde{\cal L}(\lambda).$ Then nonlocal boundary value problem~(\ref{eqP}), (\ref{eqB}) has a unique solution
 $u\in H_a^{l+2m}(K)$ for any right-hand side $\{f,\ g_{\sigma\mu}\}\in H_a^{l}(K,\ \gamma)$. This solution is represented
 in the form
 $$
  u(\omega,\ r)=(2\pi)^{-1/2}\int\limits_{-\infty+ih}^{+\infty+ih}r^{i\lambda}\tilde{\cal L}^{-1}(\lambda)
  \{\tilde F(\omega,\ \lambda),\ \tilde G_{\sigma\mu}(\lambda)\}\,d\lambda.
 $$
 Here $h=a+1-l-2m$, $\tilde F(\omega,\ \lambda)$ and $\tilde G_{\sigma\mu}(\lambda)$ are the Mellin transformations
 of the functions~$r^{2m} f(\omega,\ r)$ and~$r^{m_{\sigma\mu}}g_{\sigma\mu}(r)$ correspondingly.
\end{theorem}

Before we formulate a theorem concerning the asymptotic behavior of solutions for problem~(\ref{eqP}), (\ref{eqB}), let us
prove two Lemmas that describe solutions of the homogeneous problem.
\begin{lemma}\label{lPowerSol}
 The function
 \begin{equation}\label{eqPowerSol}
  u(\omega,\ r)=r^{i\lambda_0}\sum\limits_{q=0}^p\frac{\displaystyle 1}{\displaystyle q!}(i\ln r)^q\varphi^{(p-q)}(\omega),
 \end{equation}
 where $\varphi^{(s)}\in W^{l+2m}(b_1,\ b_2)$, $s=0,\ \dots \varkappa-1$, is a solution of homogeneous
problem~(\ref{eqP}), (\ref{eqB}) if and only if~$\lambda_0$ is an eigenvalue of the operator~$\tilde{\cal L}(\lambda)$ and
 $\varphi^{(0)},\ \dots,\ \varphi^{(\varkappa-1)}$ is a Jordan chain corresponding to the eigenvalue~$\lambda_0$;
 $p\le\varkappa-1$.
\end{lemma}
\begin{proof}
 Omitting as above the arguments~$\omega$ and~$D_\omega$ in differential operators, write
\begin{equation}\label{eqPowerSol1}
  \begin{array}{c}
  {\cal P}(D_y)u=r^{-2m}\tilde{\cal P}(rD_r)u=r^{-2m+i\lambda_0}\tilde{\cal P}(\lambda_0+rD_r)
   \sum\limits_{q=0}^p\frac{\displaystyle 1}{\displaystyle q!}(i\ln r)^q\varphi^{(p-q)}=\\
  r^{-2m+i\lambda_0}\sum\limits_{\nu=0}^p\frac{\displaystyle 1}{\displaystyle \nu!}\partial_\lambda^\nu\tilde{\cal P}(\lambda_0)
  \sum\limits_{q=\nu}^p\frac{\displaystyle 1}{\displaystyle (q-\nu)!}(i\ln r)^{q-\nu}\varphi^{(p-q)}.
 \end{array}
\end{equation}
Similarly,
\begin{equation}\label{eqPowerSol2}
  \begin{array}{c}
  B_{\sigma\mu}(D_y)u=
   r^{-m_{\sigma\mu}+i\lambda_0}
  \sum\limits_{\nu=0}^p\frac{\displaystyle 1}{\displaystyle \nu!}\partial_\lambda^\nu\tilde B_{\sigma\mu}(\lambda_0)
  \sum\limits_{q=\nu}^p\frac{\displaystyle 1}{\displaystyle (q-\nu)!}(i\ln r)^{q-\nu}\varphi^{(p-q)}.
 \end{array}
\end{equation}
Finally, consider the expression~$(B^{\cal G}_{\sigma\mu}(D_y)u)({\cal G}y)$.
\begin{equation}\label{eqPowerSol3}
  \begin{array}{c}
  (B^{\cal G}_{\sigma\mu}(D_y)u)({\cal G}y)=\\
   =r^{-m_{\sigma\mu}+i\lambda_0}\beta_\sigma^{-m_{\sigma\mu}+i\lambda_0}
  \sum\limits_{s=0}^p\frac{\displaystyle 1}{\displaystyle s!}\partial_\lambda^s\tilde B^{\cal G}_{\sigma\mu}(\lambda_0)
  \sum\limits_{q=s}^p\frac{\displaystyle 1}{\displaystyle (q-s)!}(i\ln r+i\ln\beta_{\sigma})^{q-s}
  \varphi^{(p-q)}(\omega+\omega_\sigma).
 \end{array}
\end{equation}
Applying the binomial formula to~$(i\ln r+i\ln\beta_{\sigma})^{q-s}$ and using the relation
$$
 \beta_\sigma^{-m_{\sigma\mu}+i\lambda_0}\sum\limits_{s=0}^\nu
 \frac{\displaystyle 1}{\displaystyle s!(\nu-s)!}
 \partial_\lambda^s B^{\cal G}_{\sigma\mu}(\lambda_0)(i\ln\beta_\sigma)^{\nu-s}=
\frac{\displaystyle 1}{\displaystyle \nu!}
 \partial_\lambda^\nu(\beta_\sigma^{-m_{\sigma\mu}+i\lambda}\tilde B^{\cal G}_{\sigma\mu}(\lambda))|_{\lambda=\lambda_0},
$$
we obtain from~(\ref{eqPowerSol3})
\begin{equation}\label{eqPowerSol4}
  \begin{array}{c}
  (B^{\cal G}_{\sigma\mu}(D_y)u)({\cal G}y)=\\
   =r^{-m_{\sigma\mu}+i\lambda_0}
  \sum\limits_{\nu=0}^p\frac{\displaystyle 1}{\displaystyle \nu!}\partial_\lambda^\nu
  (\beta_\sigma^{-m_{\sigma\mu}+i\lambda}\tilde B^{\cal G}_{\sigma\mu}(\lambda))|_{\lambda=\lambda_0}
  \sum\limits_{q=\nu}^p\frac{\displaystyle 1}{\displaystyle (q-\nu)!}(i\ln r)^{q-\nu}
  \varphi^{(p-q)}(\omega+\omega_\sigma).
 \end{array}
\end{equation}
Combining the summands at the same powers of~$i\ln r$ in~(\ref{eqPowerSol1}), (\ref{eqPowerSol2}), (\ref{eqPowerSol4}),
we see that the function~$u$ satisfies homogeneous problem~(\ref{eqP}), (\ref{eqB}) if and only if
$$
 \sum\limits_{h=0}^k\frac{\displaystyle 1}{\displaystyle h!}\partial_\lambda^h \tilde{\cal L}(\lambda_0)\varphi^{(k-h)}=0,\quad
 k=0,\ \dots,\ p.
$$
\end{proof}

Any solution of form~(\ref{eqPowerSol}) for homogeneous problem~(\ref{eqP}), (\ref{eqB}) is called {\it a power solution
of order~$p$} corresponding to the eigenvalue~$\lambda_0$.

Repeating the proof of Lemma~1.2~\cite{MP}, from Lemma~\ref{lPowerSol} of the present work, we derive the following statement.
\begin{lemma}\label{lBasisPowerSol}
 Let $\{\varphi^{(0,j)},\ \dots,\ \varphi^{(\varkappa_j-1,j)}:j=1,\ \dots,\ J\}$ be a canonical system of Jordan chains of
 the operator~$\tilde{\cal L}(\lambda)$ corresponding to an eigenvalue~$\lambda_0$.
 Then the functions
 \begin{equation}\label{eqPowerSolLambda}
   u^{(k,j)}(\omega,\ r)=r^{i\lambda_0}\sum\limits_{q=0}^k\frac{\displaystyle 1}{\displaystyle q!}(i\ln r)^q\varphi^{(k-q,j)}(\omega),\quad
   k=0,\ \dots,\ \varkappa_j-1,\ j=1,\ \dots,\ J,
 \end{equation}
 form a basis for the space of power solutions to homogeneous problem~(\ref{eqP}), (\ref{eqB}) corresponding to the
 eigenvalue~$\lambda_0$.
\end{lemma}

Similarly to Theorem~1.2~\cite{MP}, using Theorem~\ref{thSolv} and Lemma~\ref{lBasisPowerSol} of this work, one can prove
the following statement concerning the asymptotic representation of solutions for nonlocal problem~(\ref{eqP}), (\ref{eqB}).
\begin{theorem}\label{thAsymp}
 Let we have $\{f,\ g_{\sigma\mu}\}\in H_a^l(K,\ \gamma)\cap H_{a_1}^l(K,\ \gamma)$, where $a>a_1$. Suppose the lines
 $\Im\lambda=a_1+1-l-2m$, $\Im\lambda=a+1-l-2m$ contain no eigenvalues of the operator~$\tilde{\cal L}(\lambda)$.
 If $u$ is a solution for problem~(\ref{eqP}), (\ref{eqB}) from the space~$H_a^{l+2m}(K)$, then
 \begin{equation}\label{eqAsymp}
  u(\omega,\ r)=\sum\limits_{n=1}^{N}\sum\limits_{j=1}^{J_n}\sum\limits_{k=0}^{\varkappa_{j,n}-1}
  c_n^{(k,j)}u_n^{(k,j)}(\omega,\ r)+u_1(\omega,\ r).
 \end{equation}
 Here $\lambda_1,\ \dots,\ \lambda_N$ are eigenvalues of $\tilde{\cal L}(\lambda)$ located in the strip
 $a_1+1-l-2m<\Im\lambda<a+1-l-2m$;
 \begin{equation}\label{eqPowerSolLambda_n}
   u_n^{(k,j)}(\omega,\ r)=r^{i\lambda_n}\sum\limits_{q=0}^k
   \frac{\displaystyle 1}{\displaystyle q!}(i\ln r)^q\varphi_n^{(k-q,j)}(\omega)
 \end{equation}
 are power solutions (of order~$k$) for homogeneous problem~(\ref{eqP}), (\ref{eqB});
$$
 \{\varphi_n^{(0,j)},\ \dots,\ \varphi_n^{(\varkappa_{j,n}-1,j)}: j=1,\ \dots,\ J_n\}
$$
 is a canonical system of Jordan chains of the operator~$\tilde{\cal L}(\lambda)$ corresponding to the
 eigenvalue~$\lambda_n$, $n=1,\ \dots,\ N$; $c_n^{(k,j)}$ are some constants; $u_1$ is a solution for
 problem~(\ref{eqP}), (\ref{eqB}) from the space~$H_{a_1}^{l+2m}(K)$.
\end{theorem}

\begin{remark}\label{rAsymp}
 One can show that the formula~(\ref{eqAsymp}) is valid even if the
 line~$\Im\lambda=a+1-l-2m$ contains eigenvalues of the operator $\tilde{\cal L}(\lambda)$.
 We demand that the line~$\Im\lambda=a+1-l-2m$ has no eigenvalues, since this condition will be also
 used for studying asymptotics of solutions for the adjoint problem (Theorem~\ref{thAsympL*}).
\end{remark}

\begin{remark}
 If the conditions of Theorem~\ref{thAsymp} are fulfilled and the strip~$a_1+1-l-2m\le \Im\lambda<a+1-l-2m$
 contains no eigenvalues of the operator~$\tilde{\cal L}(\lambda)$, then the solution
 $u$ from Theorem~\ref{thAsymp} belongs to the space $H_{a_1}^{l+2m}(K)$.
\end{remark}

\section{Adjoint nonlocal problems in angles}\label{sectAdj}
{\bf 1.}
In order to calculate the coefficients~$c_\nu^{(k,j)}$ in asymptotic formula~(\ref{eqAsymp}), 
we shall need the operators that are adjoint to the operators of nonlocal problems. 

Denote $W^l[b_1,\ b_2]^*=W^l(b_1,\ b_2)^*\times{\mathbb C}^{2m}$. Consider the operator
$
 \tilde{\cal L}^*(\lambda): W^{l}[b_1,\ b_2]^*\to W^{l+2m}(b_1,\ b_2)^*,
$
which is adjoint to the operator~$\tilde{\cal L}(\bar\lambda)$ with regard to the extension of inner product 
in~$L_2(b_1,\ b_2)\times{\mathbb C}^{2m}$. The operator
$\tilde{\cal L}^*(\lambda)$ takes~$\{\psi,\ \chi_{\sigma\mu}\}\in W^{l}[b_1,\ b_2]^*$ to 
$\tilde{\cal L}^*(\lambda)\{\psi,\ \chi_{\sigma\mu}\}$ by the rule
$$ 
 <\varphi,\ \tilde{\cal L}^*(\lambda)\{\psi,\ \chi_{\sigma\mu}\}>=<\tilde{\cal P}(\bar\lambda)\varphi,\ \psi>+
\sum\limits_{\sigma=1,2}\sum\limits_{\mu=1}^m\tilde{\cal B}_{\sigma\mu}(\bar\lambda)\varphi\cdot\overline{\chi_{\sigma\mu}}\ 
\mbox{for all } \varphi\in W^{l+2m}(b_1,\ b_2).
$$
 Here and further~$<\cdot,\ \cdot>$ is a sesquilinear form on the corresponding couple of adjoint spaces.

First of all we give a remark analogous to Remark~\ref{rSmoothEigVecL}.
\begin{remark}\label{rSmoothEigVecL*} 
An eigenvector and associated vectors
$\{\psi^{(0)},\ \chi_{\sigma\mu}^{(0)}\},\ \dots,\ \{\psi^{(\varkappa-1)},\ \chi_{\sigma\mu}^{(\varkappa-1)}\}$ 
of the operator $\tilde{\cal L}^*(\lambda)$ corresponding to an eigenvalue~$\bar\lambda_0$ satisfy the equalities
\begin{equation}\label{eqDefEigVecL*}
 \sum\limits_{q=0}^\nu\frac{\displaystyle 1}{\displaystyle q!}\partial_\lambda^q \tilde{\cal L}^*(\bar\lambda_0)\
  \{\psi^{(\nu-q)},\ \chi_{\sigma\mu}^{(\nu-q)}\}=0,\quad
 \nu=0,\ \dots,\ \varkappa-1.
\end{equation}

 From equalities~(\ref{eqDefEigVecL*}) and Lemma~\ref{lSmoothL*Lambda}, it follows that the components
 $\psi^{(0)},\ \dots,\ \psi^{(\varkappa-1)}$ of an eigenvector and associated vectors of the operator~$\tilde{\cal L}^*(\lambda)$ 
 are infinitely differentiable functions in the intervals~$[b_1,\ b]$ and~$[b,\ b_2]$.
\end{remark}

Denote $H_a^l(K,\ \gamma)^*=
H_a^l(K)^*\times\prod\limits_{\sigma=1,2}\prod\limits_{\mu=1}^m H_a^{l+2m-m_{\sigma\mu}-1/2}(\gamma_\sigma)^*$.
Let ${\cal L}^*: H_a^l(K,\ \gamma)^*\to H_{a}^{l+2m}(K)^*$
be the operator adjoint to the operator~${\cal L}$ with regard to the extension of inner product 
in~$L_2(K)\times\prod\limits_{\sigma=1,2}\prod\limits_{\mu=1}^m L_2(\gamma_\sigma)$. The operator~${\cal L}^*$ 
takes~$\{v,\ w_{\sigma\mu}\}\in H_a^l(K,\ \gamma)^*$ to $\tilde{\cal L}^*\{v,\ w_{\sigma\mu}\}$ by the rule
\begin{equation}\label{eqL*Def}
 <u,\ {\cal L}^*\{v,\ w_{\sigma\mu}\}>=<{\cal P}(D_y)u,\ v>+
\sum\limits_{\sigma=1,2}\sum\limits_{\mu=1}^m<{\cal B}_{\sigma\mu}(D_y)u,\ w_{\sigma\mu}>\ \mbox{for all } u\in H_{a}^{l+2m}(K). 
\end{equation}

Consider the homogeneous equation
\begin{equation}\label{eqL*0}
 {\cal L}^*\{v,\ w_{\sigma\mu}\}=0.
\end{equation}

\begin{lemma}\label{lPowerSolL*}
  The function
 \begin{equation}\label{eqPowerSolL*} 
  \{v,\ w_{\sigma\mu}\}=\Big\{r^{i\bar\lambda_0+2m-2}\sum\limits_{q=0}^p\frac{\displaystyle 1}{\displaystyle q!}(i\ln r)^q\psi^{(p-q)},\
 r^{i\bar\lambda_0+m_{\sigma\mu}-1}\sum\limits_{q=0}^p\frac{\displaystyle 1}{\displaystyle q!}(i\ln r)^q\chi_{\sigma\mu}^{(p-q)}\Big\},
 \end{equation}
 where $\{\psi^{(s)},\ \chi^{(s)}\}\in W^{l}[b_1,\ b_2]^*$, $s=0,\ \dots \varkappa-1$, is a solution for homogeneous 
equation~(\ref{eqL*0}) if and only if~$\bar\lambda_0$ is an eigenvalue of the operator~$\tilde{\cal L}^*(\lambda)$ and
 $\{\psi^{(0)},\ \chi_{\sigma\mu}^{(0)}\},\ \dots,\ \{\psi^{(\varkappa-1)},\ \chi_{\sigma\mu}^{(\varkappa-1)}\}$ is a 
 Jordan chain corresponding to the eigenvalue~$\bar\lambda_0$; $p\le \varkappa-1$.
\end{lemma} 
\begin{proof}
By Remark~\ref{rSmoothEigVecL*} the functions~$\psi^{(s)}$, $s=0,\ \dots \varkappa-1$ belong to~$L_2(b_1,\ b_2)$. 
Therefore, for any~$u\in C_0^\infty(\bar K\backslash \{0\})$, the following identity holds:
\begin{equation}\label{eqPowerSolL*1}
 \begin{array}{c}
 <u,\ {\cal L}^* \{v,\ w_{\sigma\mu}\}>=\\
=\int\limits_{b_1}^{b_2}\int\limits_0^\infty r^{-1}\tilde{\cal P}(rD_r)u\cdot \overline{
r^{i\bar\lambda_0}\sum\limits_{q=0}^p\frac{\displaystyle 1}{\displaystyle q!}(i\ln r)^q\psi^{(p-q)}}\,dr\,d\omega+\\
+\int\limits_0^\infty \sum\limits_{\sigma=1,2}\sum\limits_{\mu=1}^m 
r^{-1} \tilde B_{\sigma\mu}(rD_r)u|_{\omega=b_\sigma}\cdot \overline{
r^{i\bar\lambda_0}\sum\limits_{q=0}^p\frac{\displaystyle 1}{\displaystyle q!}(i\ln r)^q\chi_{\sigma\mu}^{(p-q)}}\,dr+\\
+\int\limits_0^\infty \sum\limits_{\sigma=1,2}\sum\limits_{\mu=1}^m 
r^{-1} \tilde B^{\cal G}_{\sigma\mu}(rD_r)u|_{\omega=b}\cdot
 \overline{
r^{i\bar\lambda_0}\beta_\sigma^{-m_{\sigma\mu}-i\bar\lambda_0}\sum\limits_{q=0}^p\frac{\displaystyle 1}{\displaystyle q!}
(i \ln\beta_{\sigma}^{-1}+i\ln r)^q\chi_{\sigma\mu}^{(p-q)}}\,dr
 \end{array}
\end{equation}
(if we put~$r'=r\beta_\sigma^{-1}$ in the last integral, then we obtain exactly formula~(\ref{eqL*Def})).

Denote by~$\delta_{b'}=\delta_{b'}(\omega)$ the delta--function with support at the point~$b'$ ($b_1\le b'\le b_2$).
Let~$\tilde{\cal P}^*(\lambda)$, $\tilde B_{\sigma\mu}^*(\lambda)$, and~$(\tilde B^{\cal G}_{\sigma\mu})^*(\lambda)$ 
be the operators formally adjoint to~$\tilde{\cal P}(\bar\lambda)$, $\tilde B_{\sigma\mu}(\bar\lambda)$, 
and $\tilde B^{\cal G}_{\sigma\mu}(\bar\lambda)$ correspondingly.

Notice that identities of the form
$$
  \int\limits_{b_1}^{b_2} D_\omega\varphi\cdot \overline{\psi^{(p-q)}}\,d\omega=<\varphi,\ D_\omega\psi^{(p-q)}>,\quad
 D_\omega \varphi|_{\omega=b'}\cdot\overline {\chi^{(p-q)}}=<\varphi,\ D_\omega(\chi^{(p-q)}\otimes\delta_{b'})>
$$
(for~$\varphi\in W^l(b_1,\ b_2)$) generate the distributions~$D_\omega\psi^{(p-q)}$ 
and~$D_\omega(\chi^{(p-q)}\otimes\delta_{b'})$ from the space~$W^l(b_1,\ b_2)^*$.
Therefore, integrating in~(\ref{eqPowerSolL*1}) by parts (for fixed~$\omega$) and using the relations
$$
\tilde{\cal P}^*(rD_r)
\Big(r^{i\bar\lambda_0}\sum\limits_{q=0}^p\frac{\displaystyle 1}{\displaystyle q!}(i\ln r)^q\psi^{(p-q)}\Big)=
r^{i\bar\lambda_0}\sum\limits_{\nu=0}^p\frac{\displaystyle 1}{\displaystyle \nu!}\partial_\lambda^\nu\tilde{\cal P}^*(\bar\lambda_0)
  \sum\limits_{q=\nu}^p\frac{\displaystyle 1}{\displaystyle (q-\nu)!}(i\ln r)^{q-\nu}\psi^{(p-q)},
$$
$$
\begin{array}{c}
\tilde B_{\sigma\mu}^*(rD_r)
\Big(r^{i\bar\lambda_0}\sum\limits_{q=0}^p\frac{\displaystyle 1}{\displaystyle q!}(i\ln r)^q\chi_{\sigma\mu}^{(p-q)}\otimes
 \delta_{b_\sigma}\Big)=\\
=r^{i\bar\lambda_0}\sum\limits_{\nu=0}^p\frac{\displaystyle 1}{\displaystyle \nu!}
\partial_\lambda^\nu\tilde B_{\sigma\mu}^*(\bar\lambda_0)
 \Big( \sum\limits_{q=\nu}^p\frac{\displaystyle 1}{\displaystyle (q-\nu)!}(i\ln r)^{q-\nu}\chi_{\sigma\mu}^{(p-q)}\otimes
 \delta_{b_\sigma}\Big),
\end{array}
$$
$$
\begin{array}{c}
(\tilde B^{\cal G}_{\sigma\mu})^*(rD_r)
\Big(r^{i\bar\lambda_0}\beta_{\sigma}^{-m_{\sigma\mu}-i\bar\lambda_0}
\sum\limits_{q=0}^p\frac{\displaystyle 1}{\displaystyle q!}(i\ln \beta_\sigma^{-1}+i\ln r)^q\chi_{\sigma\mu}^{(p-q)}\otimes
 \delta_{b}\Big)=\\
=r^{i\bar\lambda_0}\sum\limits_{\nu=0}^p\frac{\displaystyle 1}{\displaystyle \nu!}
\partial_\lambda^\nu (\beta_\sigma^{-m_{\sigma\mu}-i\lambda}
  (\tilde B^{\cal G}_{\sigma\mu})^*(\lambda))|_{\lambda=\bar\lambda_0}
 \Big( \sum\limits_{q=\nu}^p\frac{\displaystyle 1}{\displaystyle (q-\nu)!}(i\ln r)^{q-\nu}\chi_{\sigma\mu}^{(p-q)}\otimes
 \delta_{b}\Big)
\end{array}
$$
(which are proved similarly to equalities~(\ref{eqPowerSol1}), 
(\ref{eqPowerSol2}), (\ref{eqPowerSol4})), we conclude that the function~$\{v,\ w_{\sigma\mu}\}$ 
satisfies homogeneous equation~(\ref{eqL*0}) if and only if
$$
 \sum\limits_{h=0}^k\frac{\displaystyle 1}{\displaystyle h!}\partial_\lambda^h \tilde{\cal L}^*(\bar\lambda_0)
 \{\psi^{(k-h)},\ \chi_{\sigma\mu}^{(k-h)}\}=0,\quad
 k=0,\ \dots,\ p
$$
(cf. the proof of Lemma~\ref{lPowerSol}).
\end{proof}

Any solution of form~(\ref{eqPowerSolL*}) for homogeneous equation~(\ref{eqL*0}) is called {\it a power solution of 
order~$p$} corresponding to the eigenvalue~$\bar\lambda_0$.

{\bf 2.}
Further we need a special choice of Jordan chains satisfying the conditions of biorthogonality and normalization.
Such chains are described in the following lemma.
\begin{lemma}\label{lBiortChains} 
 Suppose a canonical system of Jordan chains
 $$ 
  \{\varphi^{(0,j)},\ \dots,\ \varphi^{(\varkappa_j-1,j)}: j=1,\ \dots,\ J\}
 $$
  corresponds to an eigenvalue~$\lambda_0$ of the operator~$\tilde{\cal L}(\lambda)$.
 Then there exists a canonical system of Jordan chains
 $$ 
  \Big\{ \{\psi^{(0,j)},\ \chi_{\sigma\mu}^{(0,j)}\},\ \dots,\ \{\psi^{(\varkappa_j-1,j)},\ \chi_{\sigma\mu}^{(\varkappa_j-1,j)}\}:
  j=1,\ \dots,\ J \Big\}
 $$
 of the operator~$\tilde{\cal L}^*(\lambda)$ corresponding to the eigenvalue~$\bar\lambda_0$ such that the following relations
 hold:
 \begin{equation}\label{eqBiortChains} 
  \begin{array}{c}
  \sum\limits_{p=0}^\nu\sum\limits_{q=0}^k\frac{\displaystyle 1}{\displaystyle (\nu+k+1-p-q)!}\Big\{
  (\partial_\lambda^{\nu+k+1-p-q}\tilde{\cal P}(\lambda_0)\varphi^{(q,\xi)},\ 
  \psi^{(p,\zeta)})_{L_2(b_1,\ b_2)}+\\
   +\sum\limits_{\sigma=1,2}\sum\limits_{\mu=1}^m
  (\partial_\lambda^{\nu+k+1-p-q}\tilde{\cal B}_{\sigma\mu}(\lambda_0)\varphi^{(q,\xi)},\ 
  \chi_{\sigma\mu}^{(p,\zeta)})_{\mathbb C}\Big\}
  =
  \delta_{\xi,\zeta}\delta_{\varkappa_\xi-k-1,\nu}.
 \end{array}
 \end{equation}
 Here $\zeta,\ \xi=1,\ \dots,\ J$; $\nu=0,\ \dots,\ \varkappa_\zeta-1$; $k=0,\ \dots,\ \varkappa_\xi-1$; $\delta_{p,q}$ is the Kronecker symbol.
\end{lemma}
\begin{proof} 
 By Lemma~2.1~\cite{SkDu90}, $\lambda_0$ is a normal eigenvalue of the operator~$\tilde{\cal L}(\lambda)$, 
that is~$\dim\ker\tilde{\cal L}(\lambda_0)<\infty$, $\codim{\cal R}(\tilde{\cal L}(\lambda_0))<\infty$, and all points
of  the deleted neighborhood~$0<|\lambda-\lambda_0|<\rho$ (for sufficiently small $\rho$)
are regular ones for~$\tilde{\cal L}(\lambda)$. Thus, the necessary result
follows from Lemma~2.1~\cite{MP}.
\end{proof}

\section{Calculation of the coefficients in the asymptotics of solutions for nonlocal problems in angles}\label{sectCoefficients}
{\bf 1.} 
In this section we obtain explicit formulas for the coefficients~$c_n^{(k,j)}$ in asymptotic formula~(\ref{eqAsymp}). 
First we shall calculate the coefficients with the help of power solutions~$\{v,\ w_{\sigma\mu}\}$ for homogeneous 
equation~(\ref{eqL*0}), and then we shall obtain a representation of the coefficients in terms of the Green formula.

Let $\bar\lambda_n$ be an eigenvalue of the operator~$\tilde{\cal L}^*(\lambda)$, and let
 $$ 
  \Big\{ \{\psi_n^{(0,j)},\ \chi_{\sigma\mu,n}^{(0,j)}\},\ \dots,\ 
  \{\psi_n^{(\varkappa_{j,n}-1,j)},\ \chi_{\sigma\mu,n}^{(\varkappa_{j,n}-1,j)}\}:
  j=1,\ \dots,\ J_n \Big\}
 $$
 be Jordan chains of~$\tilde{\cal L}^*(\lambda)$ corresponding to the eigenvalue~$\bar\lambda_n$ and forming
 a canonical system.
 Consider the power solutions (of order~$\nu$) for equation~(\ref{eqL*0})
 \begin{equation}\label{eqPowerSolL*Lambda_n}
   \{v_n^{(\nu,j)},\ w_{\sigma\mu,n}^{(\nu,j)}\}=
  \Big\{r^{i\bar\lambda_n+2m-2}\sum\limits_{q=0}^\nu\frac{\displaystyle 1}{\displaystyle q!}(i\ln r)^q\psi_n^{(\nu-q,j)},\
 r^{i\bar\lambda_n+m_{\sigma\mu}-1}\sum\limits_{q=0}^\nu\frac{\displaystyle 1}{\displaystyle q!}(i\ln r)^q
 \chi_{\sigma\mu,n}^{(\nu-q,j)}\Big\},
 \end{equation}
 where $\nu=0,\ \dots,\ \varkappa_{j,n}-1$.

\begin{theorem}\label{thCoefL*}
 Let the conditions of Theorem~\ref{thAsymp} hold; then the coefficients~$c_n^{(k,j)}$ from~(\ref{eqAsymp}) are calculated
 by the formulas
 \begin{equation}\label{eqCoefL*}
  c_n^{(k,j)}=\Big(f,\ iv_n^{(\varkappa_{j,n}-k-1,j)}\Big)_{L_2(K)}+\sum\limits_{\sigma=1,2}\sum\limits_{\mu=1}^m
   \Big(g_{\sigma\mu},\ iw_{\sigma\mu,n}^{(\varkappa_{j,n}-k-1,j)}\Big)_{L_2(\gamma_\sigma)},
 \end{equation}
 where
 $\{v_n^{(\nu,j)},\ w_{\sigma\mu,n}^{(\nu,j)}\}$ is the vector defined by equality~(\ref{eqPowerSolL*Lambda_n}), and the 
 Jordan chains
 $$
 \{\varphi_n^{(0,j)},\ \dots,\ \varphi_n^{(\varkappa_{j,n}-1,j)}: j=1,\ \dots,\ J_n\},
 $$
 $$ 
  \Big\{ \{\psi_n^{(0,j)},\ \chi_{\sigma\mu,n}^{(0,j)}\},\ \dots,\ 
  \{\psi_n^{(\varkappa_{j,n}-1,j)},\ \chi_{\sigma\mu,n}^{(\varkappa_{j,n}-1,j)}\}:
  j=1,\ \dots,\ J_n \Big\},
 $$
 appearing in~(\ref{eqPowerSolLambda_n}) and~(\ref{eqPowerSolL*Lambda_n}) satisfy conditions~(\ref{eqBiortChains}) of
 biorthogonality and normalization.
\end{theorem}
Theorem~\ref{thCoefL*} is proved similar to Theorem~3.1~\cite{MP}.

\begin{remark}\label{rCoefEstimate}
 By Remark~\ref{rSmoothEigVecL*}, the functions~$\psi_n^{(\nu,j)}$ belong to the space~$L_2(b_1,\ b_2)$.
 From this and from equalities~(\ref{eqPowerSolL*Lambda_n}) and~(\ref{eqCoefL*}), it follows that 
 $$ 
  |c_n^{(k,j)}|\le c (\|\{f,\ g_{\sigma\mu}\}\|_{H_a^l(K,\ \gamma)}+\|\{f,\ g_{\sigma\mu}\}\|_{H_{a_1}^l(K,\ \gamma)})
 $$
 if~$\{f,\ g_{\sigma\mu}\}\in H_a^l(K,\ \gamma)\cap H_{a_1}^l(K,\ \gamma)$ and~$a_1+1-l-2m<\Im\lambda_n<a+1-l-2m$.
\end{remark}

From Theorems~\ref{thAsymp}, \ref{thCoefL*} and the duality conception, one can obtain the following result concerning
the asymptotics of solutions for the adjoint  problem
\begin{equation}\label{eqL*}
 {\cal L}^*\{v,\ w_{\sigma\mu}\}=\Psi.
\end{equation}
\begin{theorem}\label{thAsympL*} 
 Suppose~$\Psi\in H_a^{l+2m}(K)^*\cap H_{a_1}^{l+2m}(K)^*$, where~$a>a_1$, and the lines
 $\Im\lambda=a_1+1-l-2m$, $\Im\lambda=a+1-l-2m$ contain no eigenvalues of the operator~$\tilde{\cal L}(\lambda)$.
 If~$\{v,\ w_{\sigma\mu}\}$ is a solution for problem~(\ref{eqL*}) from the space~$H_{a_1}^{l}(K,\ \gamma)^*$, then
 \begin{equation}\label{eqAsympL*}
  \{v,\ w_{\sigma\mu}\}=\sum\limits_{n=1}^{N}\sum\limits_{j=1}^{J_n}\sum\limits_{k=0}^{\varkappa_{j,n}-1}
  d_n^{(k,j)}\{v_n^{(k,j)},\ w_{\sigma\mu,n}^{(k,j)}\}+\{V,\ W_{\sigma\mu}\}.
 \end{equation}
 Here~$\lambda_1,\ \dots,\ \lambda_N$ are eigenvalues of the operator~$\tilde{\cal L}(\lambda)$ located in the 
 strip~$a_1+1-l-2m<\Im\lambda<a+1-l-2m$; $\{v_n^{(k,j)},\ w_{\sigma\mu,n}^{(k,j)}\}$ are the vectors defined by 
 formula~(\ref{eqPowerSolL*Lambda_n});
 $d_n^{(k,j)}$ are some constants; $\{V,\ W_{\sigma\mu}\}$ is a solution for 
 problem~(\ref{eqL*}) from the space~$H_{a}^{l}(K,\ \gamma)^*$.
\end{theorem}

{\bf 2.} Consider the Green formula for nonlocal elliptic problems. For this, we introduce the set
$\gamma=\{y:\ \varphi=b\}$, which is the support of nonlocal data in problem~(\ref{eqP}), (\ref{eqB}).
Denote $K_1=\{y:\ b_1<\varphi<b\},$ $K_2=\{y:\ b<\varphi<b_2\}.$ For functions~$v(y)$ given in~$K$ we denote 
by~$v_\sigma(y)$ their restrictions on~$K_\sigma$, $\sigma=1,\ 2$. We say that~$v$ belongs 
to~${\cal C}^\infty(\bar K\backslash\{0\})$ if~$v_\sigma$ belongs to~$C^\infty(\bar K_\sigma\backslash\{0\})$,
$\sigma=1,\ 2$.

When considering the Green formula in the angle~$K$, we shall omit for short the argument~$D_y$ in differential operators.
Denote by~${\cal P}^*$ the operator formally adjoint to~${\cal P}$. 
By virtue of Theorem~4.1~\cite{GurGiess} (see also Theorem~1~\cite{GurDAN}), there exist (not unique)
  1) a system $\{B'_{\sigma\mu}\}_{\mu=1}^m$ of normal on~$\gamma_\sigma$ operators of orders~$2m-1-m'_{\sigma\mu}$ 
   with constant coefficients such that the system~$\{B_{\sigma\mu},\ B'_{\sigma\mu}\}_{\mu=1}^m$ is a Dirichlet system
   on~$\gamma_\sigma$\footnote{
 See~\cite[Chapter~2, \S 2.2]{LM} for the definition of a Dirichlet system.} of order~$2m$;
  2) a Dirichlet system~$\{B_\mu,\ B'_\mu\}_{\mu=1}^m$ on~$\gamma$ of order~$2m$ such that the 
  operators~$B_\mu$ and $B'_\mu$ are of orders~$2m-\mu$ and $m-\mu$ correspondingly.

 Whenever the choice has been made, there exist differential operators~$C_{\sigma\mu},$ $C'_{\sigma\mu},$ $T_{\nu}$, and
 $T^{\cal G}_{\sigma\nu}$  ($\sigma=1,\ 2;$ $\mu=1,\ \dots,\ m;$ $\nu=1,\ \dots,\ 2m$) with constant coefficients such that 
  I) the operators~$C_{\sigma\mu},\ C'_{\sigma\mu},\ T_{\nu}$, and $T^{\cal G}_{\sigma\nu}$ are of 
    orders~$m'_{\sigma\mu},\ 2m-1-m_{\sigma\mu},\ \nu-1$, and $\nu-1$ correspondingly;
  II) the system~$\{C_{\sigma\mu}\}_{\mu=1}^m$ covers the 
  operator~${\cal P}^*$ on~$\gamma_\sigma$ and supplements~$\{C'_{\sigma\mu}\}_{\mu=1}^m$ to a Dirichlet 
  system on~$\gamma_\sigma$ of order~$2m$;
  the system~$\{T_\nu\}_{\nu=1}^{2m}$ is a Dirichlet system on~$\gamma$ of order~$2m$;
  III) for all~$u\in C_0^\infty(\bar K\backslash\{0\}),$  $v\in {\cal C}^\infty(\bar K_\sigma\backslash\{0\})$, 
 the following Green formula is valid:
 \begin{equation}\label{eqGrP}
 \begin{array}{c}
 ({\cal P}u,\ v)_{L_2(K_\sigma)}+
 \sum\limits_{\sigma=1,2}\sum\limits_{\mu=1}^m 
 ({\cal B}_{\sigma\mu}u,\ C'_{\sigma\mu}v_\sigma|_{\gamma_\sigma})_{L_2(\gamma_\sigma)}
    +\sum\limits_{\mu=1}^m (B_{\mu}u|_{\gamma},\ {\cal T}_{\mu}v)_{L_2(\gamma)}=\\
 =\sum\limits_{\sigma=1,2} (u,\ {\cal P}^*v_\sigma)_{K_\sigma}
 +\sum\limits_{\sigma=1,2}\sum\limits_{\mu=1}^m
    (B'_{\sigma\mu}u|_{\gamma_\sigma},\ C_{\sigma\mu}v_\sigma|_{\gamma_\sigma})_{L_2(\gamma_\sigma)}
 +\sum\limits_{\mu=1}^m(B'_{\mu}u|_{\gamma},\ {\cal T}_{m+\mu}v)_{L_2(\gamma)}.
 \end{array} 
 \end{equation}
 Here
 $$
 {\cal T}_{\nu}v\equiv T_{\nu}v_1|_{\gamma}-
 T_{\nu}v_2|_{\gamma}+\sum\limits_{k=1,2} (T^{\cal G}_{k\nu}v_k) ({\cal G}^{-1}_ky)|_{\gamma},
 $$
where ${\cal G}^{-1}_{k}$ is the operator of rotation by the angle~$-\omega_k$ and expansion by~$1/\beta_k$ times in the 
plane~$\{y\}$ ($k=1,\ 2;\ \nu=1,\ \dots,\ 2m$).

Formula~(\ref{eqGrP}) generates the problem formally adjoint to problem~(\ref{eqP}), (\ref{eqB}):
\begin{equation}\label{eqP*}
  {\cal P}^*(D_y)v_\sigma=f_\sigma(y) \quad (y\in K_\sigma;\ \sigma=1,\ 2),
\end{equation}
\begin{equation}\label{eqC}
  {\cal C}_{\sigma\mu}(D_y)v\equiv C_{\sigma\mu}(D_y)v_\sigma|_{\gamma_\sigma}=g_{\sigma\mu}(y) \quad 
  (y\in\gamma_\sigma;\ \sigma=1,\ 2;\ \mu=1,\ \dots,\ m),
\end{equation}
\begin{equation}\label{eqT}
  \begin{array}{c}
  {\cal T}_{\nu}(D_y)v\equiv T_{\nu}(D_y)v_1|_{\gamma}-
 T_{\nu}(D_y)v_2|_{\gamma}+\\
 +\sum\limits_{k=1,2} (T^{\cal G}_{k\nu}(D_y)v_k) ({\cal G}^{-1}_ky)|_{\gamma}=h_{\nu}(y) \quad (y\in\gamma;\ \nu=1,\ \dots,\ 2m).
  \end{array}
\end{equation}
Problem~(\ref{eqP*})--(\ref{eqT}) is called {\it a nonlocal transmission problem in the angle~$K$}~\cite{GurGiess, GurDAN}. 

For functions~$\tilde v(\omega)$ given in the interval~$(b_1,\ b_2)$, we denote by~$\tilde v_1(\omega)$ 
and~$\tilde v_2(\omega)$ their restrictions on the intervals~$(b_1,\ b)$ and~$(b,\ b_2)$ correspondingly. 
We say that~$\tilde v$ belongs to~${\cal C}^\infty([b_1,\ b_2])$ if~$\tilde v_1$ belongs to~$C^\infty([b_1,\ b])$, 
$\tilde v_2$ belongs to~$C^\infty([b,\ b_2])$.

Write all the differential operators appearing in~(\ref{eqGrP}) in polar coordinates (omitting~$\omega$ and~$D_\omega$): 
${\cal P}(D_y)=r^{-2m}\tilde{\cal P}(rD_r),$
$B_{\sigma\mu}(D_y)=r^{-m_{\sigma\mu}} \tilde B_{\sigma\mu}(rD_r)$, etc.
By Theorem~4.3 \cite{GurGiess}, the following Green formula with the parameter~$\lambda$
is valid for any functions~$\tilde u\in C^\infty([b_1,\ b_2])$, $\tilde v\in {\cal C}^\infty([b_1,\ b_2])$:
 \begin{equation}\label{eqGrPLambda}
 \begin{array}{c}
   (\tilde{\cal P}(\lambda)\tilde u,\ \tilde v)_{L_2(b_1,\ b_2)}+\sum\limits_{\sigma=1,2}\sum\limits_{\mu=1}^m
    \tilde{\cal B}_{\sigma\mu}(\lambda)\tilde u\cdot\overline{
    \tilde C'_{\sigma\mu}(\lambda')\tilde v_\sigma|_{\omega=b_{\sigma}} }
   +\sum\limits_{\mu=1}^m
  \tilde B_{\mu}(\lambda)\tilde u|_{\omega=b}\cdot\overline{
          \tilde{\cal T}_{\mu}(\lambda')\tilde v }=\\
    =(\tilde u,\ \tilde{\cal P}^*(\lambda')\tilde v_1)_{L_2(b_1,\ b)}+(\tilde u,\ \tilde{\cal P}^*(\lambda')\tilde v_2)_{L_2(b,\ b_2)}+\\
   +\sum\limits_{\sigma=1,2}\sum\limits_{\mu=1}^m
    \tilde B'_{\sigma\mu}(\lambda)\tilde u|_{\omega=b_{\sigma}}\cdot\overline{
    \tilde C_{\sigma\mu}(\lambda')\tilde v_\sigma|_{\omega=b_{\sigma}} }
  +\sum\limits_{\mu=1}^m 
  \tilde B'_{\mu}(\lambda)\tilde u|_{\omega=b}\cdot\overline{
   \tilde{\cal T}_{m+\mu}(\lambda')\tilde v }.
 \end{array} 
 \end{equation}
 Here~$\lambda'=\bar\lambda-2i(m-1)$;
 $$
 \begin{array}{c}
 \tilde{\cal T}_{\nu}(\lambda')\tilde v=
 \tilde T_{\nu}(\lambda')\tilde v_1(\omega)|_{\omega=b}
 -\tilde T_{\nu}(\lambda')\tilde v_2(\omega)|_{\omega=b}+\\
  +\sum\limits_{k=1,2} \beta_{k}^{-i\lambda'+(\nu-1)}\tilde T_{k\nu}^{\cal G}(\lambda')
            \tilde v_k(\omega-\omega_k)|_{\omega=b}.
 \end{array}
 $$

Formula~(\ref{eqGrPLambda}) generates the problem formally adjoint to problem~(\ref{eqPLambda}), (\ref{eqBLambda}):
\begin{equation}\label{eqP*Lambda}
  \tilde{\cal P}^*(\lambda)\tilde v_1(\omega)=0 \quad (\omega\in (b_1,\ b)),\
  \tilde{\cal P}^*(\lambda)\tilde v_2(\omega)=0 \quad (\omega\in (b,\ b_2)),
\end{equation}
\begin{equation}\label{eqCLambda}
  \tilde{\cal C}_{\sigma\mu}(\lambda)\tilde v(\omega)\equiv 
  \tilde C_{\sigma\mu}(\lambda)\tilde v_\sigma(\omega)|_{\omega=\omega_\sigma}=0
  \quad (\sigma=1,\ 2;\ \mu=1,\ \dots,\ m),
\end{equation}
\begin{equation}\label{eqTLambda}
  \begin{array}{c}
   \tilde{\cal T}_{\nu}(\lambda)\tilde v(\omega)\equiv
 \tilde T_{\nu}(\lambda)\tilde v_1(\omega)|_{\omega=b}
 -\tilde T_{\nu}(\lambda)\tilde v_2(\omega)|_{\omega=b}+\\
  +\sum\limits_{k=1,2} \beta_{k}^{-i\lambda+(\nu-1)}\tilde T_{k\nu}^{\cal G}(\lambda)
            \tilde v_k(\omega-\omega_k)|_{\omega=b}=0 \quad (\nu=1,\ \dots,\ 2m).
  \end{array}
\end{equation}
Problem~(\ref{eqP*Lambda})--(\ref{eqTLambda}) is called {\it a nonlocal transmission problem
on the arc~$(b_1,\ b_2)$}~\cite{GurGiess, GurDAN}.

Notice that problem~(\ref{eqP*Lambda})--(\ref{eqTLambda}) can be also derived from problem~(\ref{eqP*})--(\ref{eqT}) 
if we put in the last one~$f_\sigma=0$, $g_{\sigma\mu}=0$, $h_\nu=0$ and do formally the Mellin transformation.

The operator
$$
 \tilde {\cal M}(\lambda): W^{l+2m}(b_1,\ b)\oplus W^{l+2m}(b,\ b_2)\to (W^{l}(b_1,\ b)\oplus W^{l}(b,\ b_2))\times
  {\mathbb C}^{2m}\times {\mathbb C}^{2m},
$$
acting by the formula
$$ 
 \tilde {\cal M}(\lambda)\tilde v=\{\tilde z,\ \tilde{\cal C}_{\sigma\mu}(\lambda)v,\ \tilde{\cal T}_{\nu}(\lambda)v\}.
$$
corresponds to problem~(\ref{eqP*Lambda})--(\ref{eqTLambda}). 
Here~$\tilde z(\omega)=\tilde{\cal P}^*(\lambda)\tilde v_1(\omega)$ for $\omega\in (b_1,\ b)$,
$\tilde z(\omega)=\tilde{\cal P}^*(\lambda)\tilde v_2(\omega)$ for $\omega\in (b,\ b_2)$. Notice that we cannot 
define~$\tilde z$ by the formula~$\tilde z(\omega)=\tilde{\cal P}^*(\lambda)\tilde v(\omega)$ for $\omega\in (b_1,\ b_2)$,
since the function~$\tilde v\in W^{l+2m}(b_1,\ b)\oplus W^{l+2m}(b,\ b_2)$ may be discontinuous at the 
point~$\omega=b$.

{\bf 3.} Now we shall establish a connection between Jordan chains of the operators~$\tilde {\cal L}^*(\lambda)$ 
and~$\tilde {\cal M}(\lambda)$. 
Put
$$
 \tilde{\cal C}'_{\sigma\mu}(\lambda)\tilde v=\tilde C'_{\sigma\mu}(\lambda)\tilde v_\sigma(\omega)|_{\omega=b_\sigma}.
$$
Repeating the proof of Proposition~2.5~\cite[Chapter~1]{NP} and using Green formula~(\ref{eqGrPLambda}) and 
Remark~\ref{rSmoothEigVecL*}, we obtain the following result.
\begin{lemma}\label{lConnectL*MLambda}
 Vectors
  $\{\psi^{(0)},\ \chi_{\sigma\mu}^{(0)}\},\ \dots,\ \{\psi^{(\varkappa-1)},\ \chi_{\sigma\mu}^{(\varkappa-1)}\}$ form
  a Jordan chain of the operator~$\tilde{\cal L}^*(\lambda)$ corresponding to an eigenvalue~$\bar\lambda_0$
  if and only if the vectors~$\psi^{(0)},\ \dots,\ \psi^{(\varkappa-1)}$ form a Jordan chain of the operator~$\tilde {\cal M}(\lambda)$
  corresponding to the eigenvalue~$\bar\lambda_0-2i(m-1)$ and the vectors~$\psi^{(k)}$ and~$\chi_{\sigma\mu}^{(k)}$ are
  connected by the relation
 $$
  \chi_{\sigma\mu}^{(k)}=\sum\limits_{r=0}^k\frac{\displaystyle 1}{\displaystyle r!}
  \partial_\lambda^r \tilde {\cal C}'_{\sigma\mu}(\bar\lambda_0-2i(m-1))\psi^{(k-r)}.
 $$
\end{lemma}

Combining Lemmas~\ref{lBiortChains} and~\ref{lConnectL*MLambda}, we get the following condition of
biorthogonality and normalization of Jordan chains in terms of the Green formula.
\begin{lemma}\label{lBiortChainsGreen} 
 Suppose a canonical system
 $$ 
  \{\varphi^{(0,j)},\ \dots,\ \varphi^{(\varkappa_j-1,j)}: j=1,\ \dots,\ J\}
 $$
 corresponds to an eigenvalue~$\lambda_0$ of the operator~$\tilde{\cal L}(\lambda)$.
 Then there exist a canonical system of Jordan chains
 $$ 
 \{\psi^{(0,j)},\ \dots,\ \psi^{(\varkappa_j-1,j)}:  j=1,\ \dots,\ J \}
 $$
 of the operator~$\tilde{\cal M}(\lambda)$ corresponding to the eigenvalue~$\bar\lambda_0-2i(m-1)$ such that
 the following relations are valid:
 \begin{equation}\label{eqBiortChainsGreen} 
  \begin{array}{c}
  \sum\limits_{p=0}^\nu\sum\limits_{q=0}^k\frac{\displaystyle 1}{\displaystyle (\nu+k+1-p-q)!}\Big\{
  (\partial_\lambda^{\nu+k+1-p-q}\tilde{\cal P}(\lambda_0)\varphi^{(q,\xi)},\ 
  \psi^{(p,\zeta)})_{L_2(b_1,\ b_2)}+\\
   +\sum\limits_{\sigma=1,2}\sum\limits_{\mu=1}^m
  (\partial_\lambda^{\nu+k+1-p-q}\tilde{\cal B}_{\sigma\mu}(\lambda_0)\varphi^{(q,\xi)},\ 
  \sum\limits_{r=0}^k\frac{\displaystyle 1}{\displaystyle r!}
  \partial_\lambda^r 
    \tilde {\cal C}'_{\sigma\mu}(\bar\lambda_0-2i(m-1))\psi^{(p-r,\zeta)})_{\mathbb C}\Big\}
  =\\
  =\delta_{\xi,\zeta}\delta_{\varkappa_\xi-k-1,\nu}.
 \end{array}
 \end{equation}
\end{lemma}

Put
$$
 {\cal C}'_{\sigma\mu}(D_y) v=C'_{\sigma\mu}(D_y)v_\sigma(y)|_{\gamma_\sigma}.
$$
Let us formulate the main result on a representation of the coefficients~$c_n^{(k,j)}$ from asymptotic 
formula~(\ref{eqAsymp}) in terms of the Green formula.
\begin{theorem}\label{thCoefGreen}
 Let conditions of Theorem~\ref{thAsymp} be fulfilled. Then the coefficients~$c_n^{(k,j)}$ from~(\ref{eqAsymp}) are 
calculated by the formulas
 \begin{equation}\label{eqCoefGreen}
  c_n^{(k,j)}=\Big(f,\ iv_n^{(\varkappa_{j,n}-k-1,j)}\Big)_{L_2(K)}+\sum\limits_{\sigma=1,2}\sum\limits_{\mu=1}^m
   \Big(g_{\sigma\mu},\ i  {\cal C}'_{\sigma\mu}(D_y)
 v_n^{(\varkappa_{j,n}-k-1,j)}\Big)_{L_2(\gamma_\sigma)}.
 \end{equation}
 Here
 $v_n^{(\nu,j)}$ is a power solution for homogeneous nonlocal transmission problem~(\ref{eqP*})--(\ref{eqT}) given by
 $$ 
  v_n^{(\nu,j)}=r^{i\bar\lambda_n+2m-2}\sum\limits_{q=0}^\nu\frac{\displaystyle 1}{\displaystyle q!}(i\ln r)^q\psi_n^{(\nu-q,j)},
 $$
 where $\{\psi_n^{(0,j)},\ \dots,\ \psi_n^{(\varkappa_{j,n}-1,j)}: j=1,\ \dots,\ J_n\}$ is a canonical system
 of Jordan chains of the operator~$\tilde{\cal M}(\lambda)$ corresponding to the eigenvalue~$\bar\lambda_n-2i(m-1)$, 
 and the chains
 $
 \{\varphi_n^{(0,j)},\ \dots,\ \varphi_n^{(\varkappa_{j,n}-1,j)}: j=1,\ \dots,\ J_n\}
 $ (appearing in~(\ref{eqPowerSolLambda_n})) and
$\{\psi_n^{(0,j)},\ \dots,\ \psi_n^{(\varkappa_{j,n}-1,j)}: j=1,\ \dots,\ J_n\}$
 satisfy conditions~(\ref{eqBiortChainsGreen}) of biorthogonality and normalization.
\end{theorem}
\begin{proof} 
Similarly to the proof of Lemma~\ref{lPowerSol}, one can show that
$v_n^{(\nu,j)}$ is a solution for homogeneous problem~(\ref{eqP*})--(\ref{eqT}) if and only if
$\psi_n^{(0,j)},\ \dots,\ \psi_n^{(\varkappa_{j,n}-1,j)}$ is a Jordan chain of the operator~$\tilde{\cal M}(\lambda)$
corresponding to the eigenvalue~$\bar\lambda_n-2i(m-1)$.

Further, we have
$$
 \begin{array}{c}
{\cal C}'_{\sigma\mu}(D_y) v_{n}^{(\nu,j)}=
 r^{i\bar\lambda_n+m_{\sigma\mu}-1}\tilde {\cal C}'_{\sigma\mu}(\bar\lambda_n-2i(m-1)+rDr)
 \sum\limits_{q=0}^\nu\frac{\displaystyle 1}{\displaystyle q!}(i\ln r)^q\psi_n^{(\nu-q,j)}=\\
 = r^{i\bar\lambda_n+m_{\sigma\mu}-1}\sum\limits_{s=0}^\nu
 \frac{\displaystyle 1}{\displaystyle s!}\partial_\lambda^s\tilde {\cal C}'_{\sigma\mu}(\bar\lambda_n-2i(m-1))
  \sum\limits_{q=s}^\nu\frac{\displaystyle 1}{\displaystyle (q-s)!}(i\ln r)^{q-s}\psi_n^{(\nu-q,j)}.
 \end{array}
$$
Changing the order of summation and applying Lemma~\ref{lConnectL*MLambda}, we get
$$
  {\cal C}'_{\sigma\mu}(D_y) v_n^{(\nu,j)}=
 r^{i\bar\lambda_n+m_{\sigma\mu}-1}\sum\limits_{q=0}^\nu\frac{\displaystyle 1}{\displaystyle q!}(i\ln r)^\nu
 \chi_{\sigma\mu,n}^{(\nu-q,j)}.
$$
Now the necessary result follows from Theorem~\ref{thCoefL*} and Lemma~\ref{lBiortChainsGreen}.
\end{proof}

{\bf 3.}
In conclusion of this section we consider the asymptotics of solutions for nonlocal problems in the angle with a special
right--hand side. Put
$$
 \begin{array}{c}
  F(\omega,\ r)=\sum\limits_{q=0}^M\frac{\displaystyle 1}{\displaystyle q!}(i\ln r)^q f^{(q)}(\omega),\
  G_{\sigma\mu}(r)=\sum\limits_{q=0}^M\frac{\displaystyle 1}{\displaystyle q!}(i\ln r)^q g_{\sigma\mu}^{(q)},\\
  \{f^{(q)},\ g_{\sigma\mu}^{(q)}\}\in W^l(b_1,\ b_2)\times{\mathbb C}^{2m}.
 \end{array}
$$
Let $\Lambda$ be some complex number. If~$\Lambda$ is an eigenvalue of the operator~$\tilde{\cal L}(\lambda)$, 
then denote by~$\varkappa(\Lambda)$ the greatest of partial multiplicities of this eigenvalue;
otherwise put~$\varkappa(\Lambda)=0$.
\begin{lemma}\label{lSpecSol}
 For problem~(\ref{eqP}), (\ref{eqB}) with right--hand side~$\{r^{i\Lambda-2m}F,\ r^{i\Lambda-m_{\sigma\mu}}G_{\sigma\mu}\}$,
 there exists a solution
 \begin{equation}\label{eqSpecSol}
  u(\omega,\ r)=
  r^{i\Lambda}\sum\limits_{q=0}^{M+\varkappa(\Lambda)}\frac{\displaystyle 1}{\displaystyle q!}(i\ln r)^q u^{(q)}(\omega),
 \end{equation}
 where $u^{(q)}\in W^{l+2m}(b_1,\ b_2)$. A solution of such a form is unique if~$\varkappa(\Lambda)=0$ (that is, 
  if~$\Lambda$ is not an eigenvalue of~$\tilde{\cal L}(\lambda)$). If~$\varkappa(\Lambda)>0$, 
  solution~(\ref{eqSpecSol}) is defined accurate to an arbitrary linear combination
  of power solutions~(\ref{eqPowerSolLambda}) corresponding to the eigenvalue~$\Lambda$.
\end{lemma}
The proof is analogous to the proof of Lemma~3.1~\cite[Chapter 3]{NP}. 

\begin{remark}
 The results of Sections~\ref{sectStatement}--\ref{sectCoefficients} are generalized for the case
 of a system of equations as well as for the case of an arbitrary number of nonlocal terms with supports on
 different rays.
\end{remark}
\section{Asymptotics of solutions for local problems in~${\mathbb R}^2\backslash\{0\}$}
{\bf 1.} When investigating nonlocal elliptic problems in plane domains, one should consider solutions not in a whole
domain~$G$ but in~$G\backslash {\cal K}$, where ${\cal K}$ is a finite set of points (see~\cite{SkMs86, SkDu91}). 
And solutions may have power singularities near the set~${\cal K}$ which corresponds to some conditions of coherence.
For studying asymptotics of solutions for such problems, we need the results of 
Sections~\ref{sectStatement}--\ref{sectCoefficients} and of this Section as well.

Let~${\cal P}(D_y)$ be a homogeneous properly elliptic differential operator of order~$2m$ with constant coefficients.

Introduce the bounded operator~${\cal P}={\cal P}(D_y): H_{a}^{l+2m}({\mathbb R}^2)\to H_a^l({\mathbb R}^2)$. We shall
study the asymptotics of solutions~$u\in H_a^{l+2m}({\mathbb R}^2)$ for the equation
\begin{equation}\label{eqP0}
  {\cal P}u=f
\end{equation}
supposing that~$f\in H_a^l({\mathbb R}^2)\cap H_{a_1}^l({\mathbb R}^2)$.

Write the operator~${\cal P}(D_y)$ in polar coordinates: 
${\cal P}(D_y)=r^{-2m}\tilde{\cal P}(\omega,\ D_\omega,\ rD_r)$. Coefficients of the 
operator~$\tilde{\cal P}(\omega,\ D_\omega,\ rD_r)$ as functions of~$\omega$ belong to the set $C^\infty_{2\pi}[0,\ 2\pi]$ of
$2\pi$-periodic infinitely differentiable functions.

Introduce the bounded operator
$
 \tilde{\cal P}(\lambda)=\tilde{\cal P}(\omega,\ D_\omega,\ \lambda):
 W_{2\pi}^{l+2m}(0,\ 2\pi)\to W_{2\pi}^l(0,\ 2\pi),
$
where~$W_{2\pi}^{l}(0,\ 2\pi)$ is a completion of the set~$C^\infty_{2\pi}[0,\ 2\pi]$ in $W^{l}(0,\ 2\pi).$

From~\cite[\S1]{SkMs86}, it follows that there exists a finite--meromorphic operator--valued 
function~$\tilde{\cal P}^{-1}(\lambda)$ such that its poles coinciding with eigenvalues of~$\tilde{\cal P}(\lambda)$ 
are located (maybe except a finite number) inside
a double angle less than~$\pi$ containing the Imaginary axis. If~$\lambda$ is not a pole, 
then~$\tilde{\cal P}^{-1}(\lambda)$ is a bounded inverse operator for~$\tilde{\cal P}(\lambda).$ 
If the line~$\Im\lambda=a+1-l-2m$ contains no poles of the operator~$\tilde{\cal P}^{-1}(\lambda)$ 
(or no eigenvalues of the operator~$\tilde{\cal P}(\lambda)$ which is the same), then by~\cite[\S1]{SkMs86}
the operator~${\cal P}$ is an isomorphism.

Using the formulated results and repeating considerations of~\cite[Chapter~3]{NP}, we shall obtain most
statements of this section.

\begin{theorem}\label{thAsymp0} 
 Suppose~$f\in H_a^l({\mathbb R}^2)\cap H_{a_1}^l({\mathbb R}^2)$, where $a>a_1$, and the 
 lines~$\Im\lambda=a_1+1-l-2m$, $\Im\lambda=a+1-l-2m$ contain no eigenvalues of the operator~$\tilde{\cal P}(\lambda)$.
 If~$u$ is a solution for problem~(\ref{eqP0}) from the space~$H_a^{l+2m}({\mathbb R}^2)$, then
 \begin{equation}\label{eqAsymp0}
  u(\omega,\ r)=\sum\limits_{n=1}^{N}\sum\limits_{j=1}^{J_n}\sum\limits_{k=0}^{\varkappa_{j,n}-1}
  c_n^{(k,j)}u_n^{(k,j)}(\omega,\ r)+u_1(\omega,\ r).
 \end{equation}
 Here $\lambda_1,\ \dots,\ \lambda_N$ are eigenvalues of~$\tilde{\cal P}(\lambda)$ located in the 
 strip~$a_1+1-l-2m<\Im\lambda<a+1-l-2m$; 
 \begin{equation}\label{eqPowerSolLambda_n0}
   u_n^{(k,j)}(\omega,\ r)=r^{i\lambda_n}\sum\limits_{q=0}^k
   \frac{\displaystyle 1}{\displaystyle q!}(i\ln r)^q\varphi_n^{(k-q,j)}(\omega)
 \end{equation}
 are power (of order~$k$) solutions for homogeneous problem~(\ref{eqP0});
$$
 \{\varphi_n^{(0,j)},\ \dots,\ \varphi_n^{(\varkappa_{j,n}-1,j)}: j=1,\ \dots,\ J_n\}
$$
is a canonical system of Jordan chains of the operator~$\tilde{\cal P}(\lambda)$ corresponding to the eigenvalue~$\lambda_n$, 
$n=1,\ \dots,\ N$; $c_n^{(k,j)}$ are some constants; $u_1$ is a solution for problem~(\ref{eqP0}) from the
space~$H_{a_1}^{l+2m}({\mathbb R}^2)$.
\end{theorem}

\begin{remark}\label{rAsymp0}
 Similarly to the case of plane angles, one can show that formula~(\ref{eqAsymp0}) is valid even
 if the line~$\Im\lambda=a+1-l-2m$ contains eigenvalues of the operator~$\tilde{\cal P}(\lambda)$.
 We demand that the line~$\Im\lambda=a+1-l-2m$ has no eigenvalues, since this condition will be also
 used for studying asymptotics of solutions for the adjoint problem (Theorem~\ref{thAsympL*0}).
\end{remark}

{\bf 2.}
Further we shall obtain explicit formulas for the coefficients~$c_n^{(k,j)}$ in asymptotic formula~(\ref{eqAsymp0}). 
First we shall calculate the coefficients with the help of power solutions for homogeneous adjoint 
equation and then we shall obtain a representation of the coefficients in terms of the Green formula.

Consider the operator~${\cal P}^*: H_a^l({\mathbb R}^2)^*\to H_{a}^{l+2m}({\mathbb R}^2)^*$ adjoint to~${\cal P}$ 
with respect to the extension of inner product in~$L_2({\mathbb R}^2)$ and the operator
$
 \tilde{\cal P}^*(\lambda): W_{2\pi}^l(0,\ 2\pi)^*\to W_{2\pi}^{l+2m}(0,\ 2\pi)^*,
$
adjoint to~$\tilde{\cal P}(\bar\lambda)$
with respect to the extension of inner product in~$L_2(0,\ 2\pi)$.

Let~$\bar\lambda_n$ be an eigenvalue of the operator~$\tilde{\cal P}^*(\lambda)$. Let
 $$ 
  \{\psi_n^{(0,j)},\ \dots,\ \psi_n^{(\varkappa_{j,n}-1,j)}: j=1,\ \dots,\ J_n \}
 $$
be Jordan chains of~$\tilde{\cal P}^*(\lambda)$ corresponding to the eigenvalue~$\bar\lambda_n$
and forming a canonical system.
Using ellipticity of the operator~$\tilde{\cal P}^*(\omega,\ D_\omega,\ \lambda)$, method of ``frozen" coefficients,
expansion of the functions~$\psi_n^{(\nu,j)}$ in the Fourier series by the functions~$e^{ik\omega}/\sqrt{2\pi}$\footnote{
 Possibility of expansion of a distribution~$\psi\in W_{2\pi}^l(0,\ 2\pi)^*$ in the Fourier series by the 
 functions~$e_k(\omega)=e^{ik\omega}/\sqrt{2\pi}$ is justified by the following 
 equalities: $<u,\ \psi>=<\sum\limits_k(u,\ e_k)_{L_2(0,\ 2\pi)}e_k,\ v>=
  (u,\ \sum\limits_k v_k e_k)_{L_2(0,\ 2\pi)}$, where~$u\in W_{2\pi}^l(0,\ 2\pi)$, $v_k=<e_k,\ v>$.
 },
 and equalities of type~(\ref{eqDefEigVecL*}), one can show that~$\psi_n^{(\nu,j)}$ are $2\pi$-periodic
 infinitely differentiable functions in the interval~$[0,\ 2\pi]$.

Consider the power solution (of order~$\nu$)
 \begin{equation}\label{eqPowerSolL*Lambda_n0} 
 \begin{array}{c}
   v_n^{(\nu,j)}=
  r^{i\bar\lambda_n+2m-2}\sum\limits_{q=0}^\nu\frac{\displaystyle 1}{\displaystyle q!}(i\ln r)^q\psi_n^{(\nu-q,j)},\quad
 \nu=0,\ \dots,\ \varkappa_{j,n}-1,
 \end{array}
 \end{equation}
 for the equation~${\cal P}^*v=0$ corresponding to the eigenvalue~$\bar\lambda_n$ of the operator~$\tilde{\cal P}^*(\lambda)$.

\begin{theorem}\label{thCoefL*0}
 Let the conditions of Theorem~\ref{thAsymp0} be fulfilled. Then the coefficients~$c_n^{(k,j)}$ from~(\ref{eqAsymp0}) 
 are calculated by the formulas
 \begin{equation}\label{eqCoefL*0}
  c_n^{(k,j)}=\Big(f,\ iv_n^{(\varkappa_{j,n}-k-1,j)}\Big)_{L_2({\mathbb R}^2)},
 \end{equation}
 where
 $v_n^{(\nu,j)}$ are defined by equalities~(\ref{eqPowerSolL*Lambda_n0}); the Jordan chains
 $\{\varphi_n^{(0,j)},\ \dots,\ \varphi_n^{(\varkappa_{j,n}-1,j)}: j=1,\ \dots,\ J_n\}$ and
 $\{\psi_n^{(0,j)},\ \dots,\ \psi_n^{(\varkappa_{j,n}-1,j)}:  j=1,\ \dots,\ J_n \Big\}$
 appearing in~(\ref{eqPowerSolLambda_n0}) and~(\ref{eqPowerSolL*Lambda_n0}) satisfy conditions
 of biorthogonality and normalization analogous to~(\ref{eqBiortChains}).
\end{theorem}

\begin{remark}\label{rCoefEstimate0}
 Since the functions~$\psi_n^{(\nu,j)}$ are infinitely differentiable,
  from equations~(\ref{eqPowerSolL*Lambda_n0}) and~(\ref{eqCoefL*0}), it follows that
 $$ 
  |c_n^{(k,j)}|\le c (\|f\|_{H_a^l({\mathbb R}^2)}+\|f\|_{H_{a_1}^l({\mathbb R}^2)})
 $$
 if $f\in H_a^l({\mathbb R}^2)\cap H_{a_1}^l({\mathbb R}^2)$ and $a_1+1-l-2m<\Im\lambda_n<a+1-l-2m$.
\end{remark}

From Theorems~\ref{thAsymp0}, \ref{thCoefL*0} and the duality conception, one can get the following result concerning
the asymptotics of solutions for the adjoint problem
\begin{equation}\label{eqL*00}
 {\cal P}^*v=\Psi.
\end{equation}

\begin{theorem}\label{thAsympL*0} 
 Suppose $\Psi\in H_a^{l+2m}({\mathbb R}^2)^*\cap H_{a_1}^{l+2m}({\mathbb R}^2)^*$, where $a>a_1$, and the lines
 $\Im\lambda=a_1+1-l-2m$, $\Im\lambda=a+1-l-2m$ contain no eigenvalues of the operator~$\tilde{\cal P}(\lambda)$.
 If~$v$ is a solution for problem~(\ref{eqL*}) from the space~$H_{a_1}^{l}({\mathbb R}^2)^*$, then
 \begin{equation}\label{eqAsympL*0}
  v=\sum\limits_{n=1}^{N}\sum\limits_{j=1}^{J_n}\sum\limits_{k=0}^{\varkappa_{j,n}-1}
  d_n^{(k,j)}v_n^{(k,j)}+V.
 \end{equation}
 Here~$\lambda_1,\ \dots,\ \lambda_N$ are eigenvalues of~$\tilde{\cal P}(\lambda)$ located in the 
  strip~$a_1+1-l-2m<\Im\lambda<a+1-l-2m$; $v_n^{(k,j)}$ are the vectors given by~(\ref{eqPowerSolL*Lambda_n0});
 $d_n^{(k,j)}$ are some constants; $V$ is a solution for problem~(\ref{eqL*00}) from the space~$H_{a}^{l}({\mathbb R}^2)^*$.
\end{theorem}

{\bf 3.} Consider the Green formula for local elliptic problems in~${\mathbb R}^2\backslash\{0\}$. It is easy to see that, for any 
 functions~$u\in C_0^\infty({\mathbb R}^2\backslash\{0\}),$  $v\in C^\infty({\mathbb R}^2\backslash\{0\})$, the following
 Green formula is valid:
 \begin{equation}\label{eqGrP0}
 \begin{array}{c}
 ({\cal P}(D_y)u,\ v)_{L_2({\mathbb R}^2)}=(u,\ {\cal P}^*(D_y)v)_{L_2({\mathbb R}^2)}.
 \end{array} 
 \end{equation}

Formula~(\ref{eqGrP0}) generates the problem formally adjoint to problem~(\ref{eqP0})
\begin{equation}\label{eqP*0}
  {\cal P}^*(D_y)v=f(y) \quad (y\in {\mathbb R}^2\backslash\{0\}),
\end{equation}

Further, it is not hard to prove that, for any 
functions~$\tilde u\in C_{2\pi}^\infty[0,\ 2\pi]$, $\tilde v\in C_{2\pi}^\infty[0,\ 2\pi]$, the following Green formula
with the parameter~$\lambda$ is valid:
 \begin{equation}\label{eqGrPLambda0}
 \begin{array}{c}
   (\tilde{\cal P}(\omega,\ D_\omega,\ \lambda)\tilde u,\ \tilde v)_{L_2(0,\ 2\pi)}
    =(\tilde u,\ \tilde{\cal P}^*(\omega,\ D_\omega,\ \lambda')\tilde v)_{L_2(0,\ 2\pi)},
 \end{array} 
 \end{equation}
 where~$\lambda'=\bar\lambda-2i(m-1)$.

Formula~(\ref{eqGrPLambda0}) generates the operator
$$
 \tilde {\cal Q}(\lambda)=\tilde{\cal P}^*(\omega,\ D_\omega,\ \lambda): 
 W_{2\pi}^{l+2m}(0,\ 2\pi)\to W_{2\pi}^{l}(0,\ 2\pi).
$$

Using Green formula~(\ref{eqGrPLambda0}) and relations of type~(\ref{eqDefEigVecL*}), one can establish a connection
between Jordan chains of the operators~$\tilde {\cal P}^*(\lambda)$ and~$\tilde {\cal Q}(\lambda)$. 
\begin{lemma}\label{lConnectL*MLambda0}
 Vectors
  $\psi^{(0)},\ \dots,\ \psi^{(\varkappa-1)}$ form a Jordan chain of the operator~$\tilde{\cal P}^*(\lambda)$
 corresponding to an eigenvalue~$\bar\lambda_0$ if and only if they form a Jordan chain of the 
 operator~$\tilde {\cal Q}(\lambda)$ corresponding to the eigenvalue~$\bar\lambda_0-2i(m-1)$.
\end{lemma}

Finally, using Lemma~\ref{lConnectL*MLambda0}, we shall formulate the main result concerning a representation of
the coefficients~$c_n^{(k,j)}$ from asymptotic formula~(\ref{eqAsymp0}) in terms of the Green formula.
\begin{theorem}\label{thCoefGreen0}
 Let the conditions of Theorem~\ref{thAsymp0} be fulfilled. Then the coefficients~$c_n^{(k,j)}$ from~(\ref{eqAsymp0}) 
 are calculated by the formula
 \begin{equation}\label{eqCoefGreen0}
  c_n^{(k,j)}=\Big(f,\ iv_n^{(\varkappa_{j,n}-k-1,j)}\Big)_{L_2({\mathbb R}^2)}.
 \end{equation}
 Here 
 $v_n^{(\nu,j)}$ is a power solution for homogeneous problem~(\ref{eqP*0}) given by formula~(\ref{eqPowerSolL*Lambda_n0});
 $\{\psi_n^{(0,j)},\ \dots,\ \psi_n^{(\varkappa_{j,n}-1,j)}: j=1,\ \dots,\ J_n\}$ is a canonical system
 of Jordan chains of the operator~$\tilde{\cal Q}(\lambda)$ corresponding to the eigenvalue~$\bar\lambda_n-2i(m-1)$; the
 chains
 $
 \{\varphi_n^{(0,j)},\ \dots,\ \varphi_n^{(\varkappa_{j,n}-1,j)}: j=1,\ \dots,\ J_n\}
 $ (appearing in~(\ref{eqPowerSolLambda_n0})) and
$\{\psi_n^{(0,j)},\ \dots,\ \psi_n^{(\varkappa_{j,n}-1,j)}: j=1,\ \dots,\ J_n\}$
 satisfy the conditions of biorthogonality and normalization analogous to~(\ref{eqBiortChainsGreen}).
\end{theorem}

{\bf 4.}
When investigating asymptotics of solutions for nonlocal problems in bounded domains, we need a result
on the asymptotics of solutions for adjoint local problems in~${\mathbb R}^2\backslash\{0\}$ with a special right--hand side. 
We pay attention to the distinct from the model problem in the angle where we needed a result on the asymptotics
of solutions for the origin (but not adjoint) problem with a special right--hand side.

Let~$\Lambda$ be some complex number. If~$\bar\Lambda$ is an eigenvalue of the 
operator~$\tilde{\cal P}^*(\lambda)$, then we denote by~$\varkappa(\bar\Lambda)$ the greatest of partial multiplicities
of this eigenvalue. Otherwise we put~$\varkappa(\bar\Lambda)=0$.
\begin{lemma}\label{lSpecSol*0}
 For problem~(\ref{eqP*0}) with right--hand side
 $\Psi=r^{i\bar\Lambda-2}\sum\limits_{q=0}^M\frac{\displaystyle 1}{\displaystyle q!}(i\ln r)^q \Psi^{(q)},$
 $\Psi^{(q)}\in W_{2\pi}^{l+2m}(0,\ 2\pi)^*$,
 there exist a solution
 \begin{equation}\label{eqSpecSol*0}
  v=r^{i\bar\Lambda+2m-2}\sum\limits_{q=0}^{M+\varkappa(\bar\Lambda)}
 \frac{\displaystyle 1}{\displaystyle q!}(i\ln r)^q v^{(q)},
 \end{equation}
 where~$v^{(q)}\in W_{2\pi}^{l}(0,\ 2\pi)^*$. A solution of such a form is unique if~$\varkappa(\bar\Lambda)=0$ (that is, 
 if~$\bar\Lambda$ is not an eigenvalue of~$\tilde{\cal P}^*(\lambda)$). If~$\varkappa(\bar\Lambda)>0$, then 
 solution~(\ref{eqSpecSol*0}) is defined accurate to an arbitrary linear combination
 of power solutions~(\ref{eqPowerSolL*Lambda_n0}) corresponding to the eigenvalue~$\bar\Lambda$.
\end{lemma}
\begin{proof}
The idea of the proof is analogous to the one of the proof of Lemma~3.1~\cite[Chapter~3]{NP}. To complete the picture
we shall give a plan of the proof.
One should substitute formula~(\ref{eqSpecSol*0}) of the solution into the equation
$$
  {\cal P}^*v=r^{i\bar\Lambda-2}\sum\limits_{q=0}^M\frac{\displaystyle 1}{\displaystyle q!}(i\ln r)^q \Psi^{(q)},
$$
reduce the factor~$r^{i\bar\Lambda-2}$, and gather the coefficients at the same powers
of~$i\ln r$.
As a result, one obtains a system of~$M+\varkappa(\bar\Lambda)$ equations, from which one finds unknown
~$v^{(q)}$. The statement that a solution of form~(\ref{eqSpecSol*0}) is unique 
(for~$\varkappa(\bar\Lambda)=0$) or defined accurate to an arbitrary linear combination
 of power solutions~(\ref{eqPowerSolL*Lambda_n0}) corresponding to the eigenvalue~$\bar\Lambda$ 
(for~$\varkappa(\bar\Lambda)>0$) follows from the result analogous to Lemma~1.3~\cite[Chapter~3]{NP}, which
restricts the freedom in choosing power solutions for the equation~${\cal P}^*v=0$.
\end{proof}

\appendix

\renewcommand{\section}{\@startsection{section}{2}{0pt}{-3.5ex plus -1ex minus -.2ex}{1ex}{\bf Appendix }}

\section{Smoothness of solutions to nonlocal problems for ordinary differential equations}\label{appendSmooth}
In this Appendix, we establish two auxiliary lemmas concerning smoothness of the above-mentioned problems.
These lemmas are necessary to prove smoothness of eigenvectors and associated vectors of
nonlocal elliptic problems.

Let $\tilde{\cal P}(\lambda),$ $\tilde B_{\sigma\mu}(\lambda),$ $\tilde B^{\cal G}_{\sigma\mu}(\lambda),$
$\tilde {\cal B}_{\sigma\mu}(\lambda)$ be the differential operators
defined in Section~\ref{sectStatement}.

Consider the operator
$$
 \begin{array}{c}
 \tilde{\cal L}_{(l)}(\lambda)=\{\tilde{\cal P}(\lambda),\ 
  \tilde{\cal B}_{\sigma\mu}(\lambda)\}:W^{l+2m}(b_1,\ b_2)
  \to W^{l}[b_1,\ b_2]=W^l(b_{1},\ b_{2}) \times{\mathbb C}^{2m}.
 \end{array}
$$
We study smoothness of solutions for the nonlocal problem
\begin{equation}\label{eqLlLambda}
 \tilde{\cal L}_{(l)}(\lambda) u=\{f,\ g_{\sigma\mu}\}.
\end{equation}

\begin{lemma}\label{lSmoothLLambda}
  Let $u\in W^{l+2m}(b_1,\ b_2)$ be a solution for problem~(\ref{eqLlLambda}) with right--hand side
  $\{f,\ g_{\sigma\mu}\}\in W^{l+k}(b_1,\ b_2)$. Then $u\in W^{l+2m+k}(b_1,\ b_2)$.
\end{lemma}
\begin{proof}
The function $u(\omega)$ is a solution for the problem
$$
 \begin{array}{c}
 \tilde{\cal P}(\lambda)u(\omega)=f(\omega)\quad (\omega\in (b_1,\ b_2)),\\
 \tilde B_{\sigma\mu}(\lambda)u(\omega)|_{\omega=b_\sigma}=g_{\sigma\mu}-
 \beta_\sigma^{-m_{\sigma\mu}+i\lambda} \tilde B^{\cal G}_{\sigma\mu}(\lambda)
              u(\omega+\omega_\sigma,\ \lambda)|_{\omega=b_\sigma},\\
 \sigma=1,\ 2;\ \mu=1,\ \dots,\ m.
 \end{array}
$$
Therefore, applying Theorem~5.1~\cite[Chapter~2]{LM}, we obtain~$u\in W^{l+2m+k}(b_1,\ b_2)$.
\end{proof}

Consider the operator~$\tilde{\cal L}^*_{(l)}(\lambda):W^{l}[b_1,\ b_2]^*\to W^{l+2m}(b_1,\ b_2)^*$, adjoint to the 
operator~$\tilde{\cal L}_{(l)}(\bar\lambda)$ with regard to the extension of inner product 
in~$L_2(b_1,\ b_2)\times{\mathbb C}^{2m}$ (see Section~\ref{sectAdj}).

We shall investigate smoothness of solutions for the adjoint nonlocal problem
\begin{equation}\label{eqL*lLambda}
 \tilde{\cal L}^*_{(l)}(\lambda) \{v,\ w_{\sigma\mu}\}=\Psi.
\end{equation}

\begin{lemma}\label{lSmoothL*Lambda}
  Let $\{v,\ w_{\sigma\mu}\}\in W^{l}[b_1,\ b_2]^*$ be a solution for problem~(\ref{eqL*lLambda}) with right--hand side
  $\Psi\in \left\{
   \begin{array}{l}
    W^{2m-k}(b_1,\ b_2)^*\ \mbox{for } 0<k<2m,\\
    W^{-2m+k}(b_1,\ b)\oplus W^{-2m+k}(b,\ b_2)\ \mbox{for } k\ge 2m.\\
   \end{array}
    \right.
  $
  
Then $v\in  W^{k}(b_1,\ b)\oplus W^{k}(b,\ b_2)$.
\end{lemma}
\begin{proof}
1) First, let we have $l=0$. Denote~$\tilde{\cal L}(\lambda)=\tilde{\cal L}_{(0)}(\lambda)$,
  $\tilde{\cal L}^*(\lambda)=\tilde{\cal L}^*_{(0)}(\lambda)$.

Introduce the auxiliary operator
$\tilde{\cal L}^*_{\cal G}(\lambda):L_2(b_1,\ b_2)\times {\mathbb C}^{2m}\times {\mathbb C}^{2m}  \to W^{2m}(b_1,\ b_2)^*$
taking~$\{v,\ w_{\sigma\mu},\ w'_{\sigma\mu}\}$ to $\tilde{\cal L}^*_{\cal G}(\lambda)\{v,\ w_{\sigma\mu},\ w'_{\sigma\mu}\}$ 
by the rule
$$
 \begin{array}{c}
 <u,\ \tilde{\cal L}^*_{\cal G}(\lambda)\{v,\ w_{\sigma\mu},\ w'_{\sigma\mu}\}>=
 (\tilde{\cal P}(\lambda)u,\ v)_{L_2(b_1,\ b_2)}
 +\sum\limits_{\sigma=1,2}\sum\limits_{\mu=1}^m 
         \tilde B_{\sigma\mu}(\lambda)u|_{\omega=b_\sigma}\cdot \overline{w_{\sigma\mu}}+\\
 +\sum\limits_{\sigma=1,2}\sum\limits_{\mu=1}^m 
      \beta_\sigma^{-m_{\sigma\mu}+i\lambda} \tilde B^{\cal G}_{\sigma\mu}(\lambda)
              u|_{\omega=b}\cdot \overline{w'_{\sigma\mu}} \quad \mbox{for all } u\in W^{2m}(b_1,\ b_2).
 \end{array}
$$
Clearly,
$$
 \tilde{\cal L}^*_{\cal G}(\lambda)\{v,\ w_{\sigma\mu},\ w_{\sigma\mu}\}=
\tilde{\cal L}^*(\lambda)\{v,\ w_{\sigma\mu}\}.
$$
Introduce infinitely differentiable functions~$\zeta_\sigma(\omega)$ ($\sigma=1,\ 2$), $\zeta(\omega)$,
$$
 \begin{array}{c}
 \zeta_{\sigma}(\omega)=1\ \mbox{for } |b_\sigma-\omega|<|b_\sigma-b|/4,\ 
 \zeta_{\sigma}(\omega)=0\ \mbox{for } |b_\sigma-\omega|>|b_\sigma-b|/2;\\
 \zeta(\omega)=1-\zeta_1(\omega)-\zeta_2(\omega).
 \end{array}
$$

2) Consider the expression~$\tilde{\cal L}^*_{\cal G}(\lambda)(\zeta_1\{v,\ w_{\sigma\mu},\ w_{\sigma\mu}\})$. Then we have
$$
  \begin{array}{c}
 <u,\ \tilde{\cal L}^*_{\cal G}(\lambda)(\zeta_1\{v,\ w_{\sigma\mu},\ w_{\sigma\mu}\})>=
 (\tilde{\cal P}(\lambda)u,\ \zeta_1 v)_{L_2(b_1,\ b_2)}
 +\sum\limits_{\mu=1}^m 
         \tilde B_{1\mu}(\lambda)u|_{\omega=b_1}\cdot \overline{w_{1\mu}}\\
  \mbox{for all } u\in W^{2m}(b_1,\ b_2).
 \end{array}
$$
Besides, from Leibniz's formula, it follows that
$\tilde{\cal L}^*_{\cal G}(\lambda)(\zeta_1\{v,\ w_{\sigma\mu},\ w_{\sigma\mu}\})\in W^{2m-1}(b_1,\ b_2)^*$
since
$$\zeta_1\tilde{\cal L}^*_{\cal G}(\lambda)\{v,\ w_{\sigma\mu},\ w_{\sigma\mu}\}=
\zeta_1\tilde{\cal L}^*(\lambda)\{v,\ w_{\sigma\mu}\}\in W^{2m-1}(b_1,\ b_2)^*
$$
and $v\in L_2(b_1,\ b_2)$. Therefore we can use Theorem~5.1~\cite[Chapter~2]{LM} which yields
$\zeta_1 v\in W^1(b_1,\ b_2)$.

Similarly, we get $\zeta_2 v\in W^1(b_1,\ b_2)$.

3) Consider the expression~$\tilde{\cal L}^*_{\cal G}(\lambda)(\zeta\{v,\ w_{\sigma\mu},\ w_{\sigma\mu}\})$. Then we have
$$
  \begin{array}{c}
 <u,\ \tilde{\cal L}^*_{\cal G}(\lambda)(\zeta\{v,\ w_{\sigma\mu},\ w_{\sigma\mu}\})>=
 (\tilde{\cal P}(\lambda)u,\ \zeta v)_{L_2(-\infty,\ b)}\\
   \mbox{for all } u\in C_0^\infty(-\infty,\ b),
 \end{array}
$$
where $v(\omega)$ is extended by zero for~$\omega\le b_1$. Analogously to the above, we have
$\tilde{\cal L}^*_{\cal G}(\lambda)(\zeta\{v,\ w_{\sigma\mu},\ w_{\sigma\mu}\})\in W^{-2m+1}(-\infty,\ b)$.\footnote{
 $W^{-s}(-\infty,\ b)$, $s\ge0$, is the space adjoint
 to~$\mathaccent23 W^s(-\infty,\ b)$, where~$\mathaccent23 W^s(-\infty,\ b)$ is a completion of the 
 set~$C_0^\infty(-\infty,\ b)$ in the norm $\|u\|=\left(\sum\limits_{j=0}^s\int\limits_{-\infty}^b 
  \left|\frac{\displaystyle d^{j} v}{\displaystyle d\omega^j}\right|^2\,d\omega\right)^{1/2}$.
} 
From this,
from ellipticity of the operator~$\tilde{\cal P}(\lambda)$, and relation~$v\in L_2(-\infty,\ b)$, it follows that
the generalized derivative~$\frac{\displaystyle d^{2m} (\zeta v)}{\displaystyle d\omega^{2m}}$ belongs to the 
space~$W^{-2m+1}(-\infty,\ b)$.
Therefore, by Lemma~12.3 \cite[Chapter~1]{LM}, we have $\zeta v\in W^1(-\infty,\ b)$. Similarly, one can prove
that~$\zeta v\in W^1(b,\ +\infty)$. Combining this with item~2) of the proof, we obtain~$v\in W^{1}(b_1,\ b)\oplus W^{1}(b,\ b_2)$.

Repeating the described procedure, after a finite number of steps we shall get~$v\in W^{k}(b_1,\ b)\oplus W^{k}(b,\ b_2)$.

4) Finally, consider the case of an arbitrary~$l\ge 0$. From Lemma~\ref{lSmoothLLambda}, it follows that
\begin{equation}\label{eqSmoothL*Lambda1}
 {\cal R}(\tilde{\cal L}_{(l)}(\lambda))= {\cal R}(\tilde{\cal L}_{(0)}(\lambda))\cap W^l[b_1,\ b_2].
\end{equation}
Besides, by Lemma~2.1~\cite{SkDu90}, ${\cal R}(\tilde{\cal L}_{(l)}(\lambda))$ is closed and
$\codim{\cal R}(\tilde{\cal L}_{(l)}(\lambda))$ is finite. From this and from~(\ref{eqSmoothL*Lambda1}), it follows that
the embedding~$W^l[b_1,\ b_2]$ into $W^0[b_1,\ b_2]$ induces the isomorphism
between the coset spaces
$W^l[b_1,\ b_2]/{\cal R}(\tilde{\cal L}_{(l)}(\lambda))$ and $W^0[b_1,\ b_2]/{\cal R}(\tilde{\cal L}_{(0)}(\lambda))$. 

Thus, we have $\codim{\cal R}(\tilde{\cal L}_{(l)}(\lambda))=\codim{\cal R}(\tilde{\cal L}_{(0)}(\lambda))$, and hence
$\dim\ker(\tilde{\cal L}^*_{(l)}(\lambda))=\dim\ker(\tilde{\cal L}^*_{(0)}(\lambda))$. From this and from
the evident embedding
$\ker(\tilde{\cal L}^*_{(0)}(\lambda))\subset\ker(\tilde{\cal L}^*_{(l)}(\lambda))$, we obtain
$\ker(\tilde{\cal L}^*_{(l)}(\lambda))=\ker(\tilde{\cal L}^*_{(0)}(\lambda))$.

Further, since~$\Psi\in {\cal R}(\tilde{\cal L}^*_{(l)}(\lambda))$, we have
$$ 
 <u,\ \Psi>=0 \quad \mbox{for all } u\in \ker(\tilde{\cal L}_{(l)}(\lambda)).
$$
But from Lemma~\ref{lSmoothLLambda}, it follows that~$\ker(\tilde{\cal L}_{(l)}(\lambda))=\ker(\tilde{\cal L}_{(0)}(\lambda))$.
Therefore,
$$ 
 <u,\ \Psi>=0 \quad \mbox{for all } u\in \ker(\tilde{\cal L}_{(0)}(\lambda)).
$$
Hence, we have $\Psi\in {\cal R}(\tilde{\cal L}^*_{(0)}(\lambda))$ since $\Psi\in W^{2m}(b_1,\ b_2)^*$ by assumption. 
Let~$\{f,\ g_{\sigma\mu}\}\in W^0[b_1,\ b_2]^*=W^0[b_1,\ b_2]$ be some
solution of the problem~$\tilde{\cal L}^*_{(0)}(\lambda) \{f,\ g_{\sigma\mu}\}=\Psi$. By proved, we have
$f\in W^{k}(b_1,\ b)\oplus W^{k}(b,\ b_2)$.

Clearly, $\{f,\ g_{\sigma\mu}\}$ is also a solution of the problem~$\tilde{\cal L}^*_{(l)}(\lambda) \{f,\ g_{\sigma\mu}\}=\Psi$;
therefore,
$$
 \{v,\ w_{\sigma\mu}\}-\{f,\ g_{\sigma\mu}\}\in \ker(\tilde{\cal L}^*_{(l)}(\lambda))=
 \ker(\tilde{\cal L}^*_{(0)}(\lambda)).
$$
Hence, $v$ also belongs to~$W^{k}(b_1,\ b)\oplus W^{k}(b,\ b_2)$.
\end{proof}

The author is profoundly grateful to A.L. Skubachevskii for his attention to this work.
\renewcommand{\section}{\@startsection{section}{2}{0pt}{-3.5ex plus -1ex minus -.2ex}{1ex}{\bf}}

\makeatother
\end{document}